\documentclass[10pt]{article}
\usepackage{amssymb,amsmath}
\usepackage{color}

\usepackage[italian, english]{babel} 
\usepackage[OT2, T1]{fontenc} 
\usepackage[utf8]{inputenc} 

\usepackage{hyperref}
\hypersetup{hidelinks} 

\newtheorem{theorem}{Theorem}[section]
\newtheorem{proposition}[theorem]{Proposition}
\newtheorem{corollary}[theorem]{Corollary}
\newtheorem{lemma}[theorem]{Lemma}
\newtheorem{remark}[theorem]{Remark}
\newtheorem{definition}[theorem]{Definition}

\def\claim#1.{\noindent {\bf #1.}}

\def\flushright#1{{\unskip\nobreak\hfil\penalty50\hskip2em\hbox{}\nobreak\hfil%
#1\parfillskip=0pt\finalhyphendemerits=0\par}}

\def\bull{\vrule height 1.8ex width 1.0ex depth .1ex }
\def\QED{\ifmmode\eqno\hbox{$\bull$}\else\flushright{\hbox{$\bull$}}\fi}

\newcommand{\parag}[1]{\left\{ \begin{array}{rcl} #1 \end{array}\right.} 

\newcommand{\R}{\mathbb{R}}
\newcommand{\N}{\mathbb{N}}
\newcommand{\Z}{\mathbb{Z}}
\newcommand{\G}{{\mathbb G}}

\newcommand{\nequiv}{\not \equiv}
\newcommand{\mc}[1]{\mathcal{#1}}
\newcommand{\genus}{{\hbox{genus}}}

\def\half{{1\over 2}}

\def\abs#1{|#1|}

\def\norm#1{\|#1\|}
\def\dist{\mathop{\hbox{dist}}\nolimits}
\def\intRN{\int_{\R^N}}

\def\epsilon{\varepsilon}

\def\subseteq{\subset}

\usepackage{scalerel}
\usepackage{tikz}
\usetikzlibrary{svg.path}
\definecolor{orcidlogocol}{HTML}{A6CE39}
\tikzset{
 orcidlogo/.pic={
 \fill[orcidlogocol] svg{M256,128c0,70.7-57.3,128-128,128C57.3,256,0,198.7,0,128C0,57.3,57.3,0,128,0C198.7,0,256,57.3,256,128z};
 \fill[white] svg{M86.3,186.2H70.9V79.1h15.4v48.4V186.2z}
 svg{M108.9,79.1h41.6c39.6,0,57,28.3,57,53.6c0,27.5-21.5,53.6-56.8,53.6h-41.8V79.1z M124.3,172.4h24.5c34.9,0,42.9-26.5,42.9-39.7c0-21.5-13.7-39.7-43.7-39.7h-23.7V172.4z}
 svg{M88.7,56.8c0,5.5-4.5,10.1-10.1,10.1c-5.6,0-10.1-4.6-10.1-10.1c0-5.6,4.5-10.1,10.1-10.1C84.2,46.7,88.7,51.3,88.7,56.8z};
 }
}
\newcommand\orcidicon[1]{\href{https://orcid.org/#1}{\mbox{\scalerel*{
\begin{tikzpicture}[yscale=-1,transform shape]
\pic{orcidlogo};
\end{tikzpicture}
}{|}}}}

\begin{document}

\title {
Normalized solutions for fractional nonlinear scalar field equations via Lagrangian formulation}

\author{
Silvia Cingolani \orcidicon{0000-0002-3680-9106}\\
Dipartimento di Matematica\\
Universit\`{a} degli Studi di Bari Aldo Moro\\
Via E. Orabona 4, 70125 Bari, Italy\\
\href{mailto:silvia.cingolani@uniba.it}{silvia.cingolani@uniba.it} 
\\ 
\\
Marco Gallo \orcidicon{0000-0002-3141-9598}\\
	Dipartimento di Matematica\\
	Universit\`{a} degli Studi di Bari Aldo Moro\\
	Via E. Orabona 4, 70125 Bari, Italy\\ 
	\href{mailto:marco.gallo@uniba.it}{marco.gallo@uniba.it}
\\ 
\\
Kazunaga Tanaka \orcidicon{0000-0002-1144-1536}\\
Department of Mathematics\\
School of Science and Engineering\\
Waseda University\\
3-4-1 Ohkubo, Shijuku-ku, Tokyo 169-8555, Japan\\
\href{mailto:kazunaga@waseda.jp}{kazunaga@waseda.jp} 
}

\date{}

\maketitle


\begin{abstract}
\noindent 
We study existence of solutions for the fractional problem
\begin{equation*}
(P_m) \quad	\parag{
	(-\Delta)^{s} u + \mu u &=g(u) & \; \text{in $\mathbb{R}^N$}, \cr
	\int_{\mathbb{R}^N} u^2 dx &= m, & \cr
	 u \in H^s_r&(\mathbb{R}^N), &
	} 
\end{equation*}
where $N\geq 2$, $s\in (0,1)$, $m>0$, $\mu$ is an unknown Lagrange multiplier and $g \in C(\mathbb{R}, \mathbb{R})$ satisfies Berestycki-Lions type conditions.
Using a Lagrangian formulation of the problem $(P_m)$, we prove the existence of a weak solution with prescribed mass when $g$ has $L^2$ subcritical growth. 
The approach relies on the construction of a minimax structure, by means of a \emph{Pohozaev's mountain} in a product space and some deformation arguments under a new version of the Palais-Smale condition introduced in \cite{HT,IT}.
A multiplicity result of infinitely many normalized solutions is also obtained if $g$ is odd.%
\end{abstract}

 \medskip



\tableofcontents


\setcounter{equation}{0} 
\section{Introduction}\label{section:1}

In 1948, following a suggestion by P.A.M. Dirac, R.P. Feynman proposed a new suggestive description of the time evolution of the state of a non-relativistic quantum particle. 
According to Feynman, the wave function solution of the Schr\"odinger equation should be given by a heuristic integral over the space of paths. 
The classical notion of a single, unique classical trajectory for a system is replaced by a functional integral over an infinity of quantum-mechanically possible trajectories. 
Following Feynman's path integral approach to quantum mechanics, Laskin \cite{Las0} generalized the path integral over Brownian motions (random motion seen in swirling gas molecules) to L\'evy flights (a mix of long trajectories and short, random movements found in turbulent fluids) and derived the fractional nonlinear Schr\"odinger ((fNLS) for short) equation 
\begin{equation}\label{fNLS}
i \partial_t \psi = (-\Delta)^s \psi - g(\psi), \quad (t,x) \in \R \times \R^N ,
\end{equation}
where $\psi(t,x)$ is a complex wave, $s \in (0,1)$, the symbol $(-\Delta)^{s}$ denotes the fractional power of the Laplace operator and $g$ is a Gauge invariant nonlinearity, i.e. $g(e^{i \theta} \rho) = e^{i \theta} g(\rho)$ for any $\rho$, $\theta \in \R$.

In 2015 a first optical realization of the fractional Schr\"odinger equation, based on transverse light dynamics in aspherical optical cavities, was achieved by Longhi \cite{Lon0}. 
Subsequently, the propagation dynamics of wave packets were reported in Kerr nonlinearity, with constant or double-barrier potential. 
Numerical results showed the existence of solitons for (fNLS) equations where the L\'evy index $s$ and the saturation parameter can significantly affect the stability of these solitons \cite{KSM,WCCYZH,LMM,YL}.
Numerous other applications of the (fNLS) equation arise in the physical sciences, ranging from models of boson stars \cite{FJL} to geo-hydrology \cite{Ata0}, from charge transport in biopolymers, like DNA \cite{KLS} to anomalous diffusion phenomena \cite{BuV,Vaz0}, from water wave dynamics \cite{IP} to jump processes in probability theory with applications to financial mathematics (see also \cite{DPV} and the references therein).

\label{pag_comm_physics1}

Moreover, fractional integrals and derivatives in the calculation methods have been used for the explanation of physical phenomena which do not comply with the laws of classical statistical physics, for instance in modelling Bose-Einstein condensates. 
It is known that Bose-Einstein condensation, theoretically discovered in 1924 and observed experimentally with alkali metals (rubidium and sodium atoms) in 1995, represents a topical subject due to the explanation of quantum effects seen on a macroscopic scale, transmission of matter and the behaviour of superconductivity and superfluids. 
In this respect, not only experimental studies are important but theoretical studies too, which lead to the analysis of class of (fNLS) equations (also known as fractional Gross-Pitaievskii equations). 
Numerical simulations show existence of standing waves solutions, having a soliton behaviour and bound states \cite{DZ,ZHSWW}, including mass conservation, energy conservation and dispersion relation, in which the fractional order exponent influences the shape of the state.

From a mathematical view point, when searching for standing waves to \eqref{fNLS}, i.e. factorized solutions $\psi (t,x) = e^{i \mu t } u(x)$, $\mu \in \R$, two possible directions can be pursued. 
A first possibility is to study \eqref{fNLS} with a prescribed frequency $\mu$ and free mass. 
This approach, which we call the \emph{unconstrained} problem, has been deeply developed, for example by \cite{FQT,CW,BKS,Iko1,Iko2}. 
The literature concerning the local version of the unconstrained problem starts from the seminal paper of Berestycki-Lions \cite{BL1} and it is so large that we do not even make an attempt to summarize it.
 
A second approach is to prescribe the mass of $u$, thus conserved by $\psi$ in time
$$\int_{\R^N} |\psi(x,t)|^2 \, dx= m, \quad \forall \, t\in \R$$
and let the frequency $\mu$ to be free, becoming an unknown. 
This second approach is of considerable significance in physics, not only for the information on the mass itself, but also because the mass may also have specific meaning, such as the power supply in nonlinear optics, or the total number of atoms in Bose-Einstein condensation. 
Moreover, it can give better insights into the dynamical properties, such as the orbital stability or instability of solutions of \eqref{fNLS}.

 In a local framework ($s=1$) the seminal contribution to the study of \emph{constrained} problems is due to Stuart \cite{Stu0} and Cazenave and Lions \cite{CL}. 
 See \cite{Jea0,BV,Shi0,HT} for more recent contributions in the local case (see also \cite{CT} for a NLS equation with a nonlocal source term).
 
In the nonlocal case, it remains an open problem to derive analytically the existence of infinitely many bound states with higher energy, including mass conservation.

The present work is dedicated to the study of standing waves solutions of \eqref{fNLS} with prescribed mass by means of a new variational method. 
Namely, we are interested to seek for radially symmetric solutions of the fractional problem
\begin{equation*} \label{problem}
(P_m) 
\quad \parag{
(-\Delta)^{s} u + \mu u &=& g(u) \quad \hbox{in $\R^N$,} \cr
\int_{\R^N} |u|^2 dx &=& m, 
}
\end{equation*}
where $N\geq 2$, $s\in (0,1)$, $m>0$ and $\mu$ is a Lagrange multiplier. 
We assume that the function $g$ satisfies the following Berestycki-Lions type conditions:
\vskip2pt
\begin{itemize}
	\item[(g1)] $g : \R \to \R$ continuous and $\lim_{t \to 0} \frac{g(t)}{t}=0$,
	\item[(g2)] $\lim_{|t| \to \infty} \frac{g(t)}{|t|^p} =0$ where $p = 1 + \frac{4s}{N}$,
	\item[(g3)] there exists $t_0>0$ such that $G(t_0) >0$,
\end{itemize}
\vskip 2pt
where $G(t)= \int_0^t g(\tau) d\tau$. 
We remark that the exponent $p = 1 + \frac{4s}{N}$ appears as a $L^2$ critical exponent for the nonlinear fractional equations with $L^2$ constraint and thus assumption (g2) means that $g$ has $L^2$ subcritical growth.

Some nonlinear models satisfying (g1)--(g3) are given by pure powers $g(t)=|t|^{q-2}t$, with $q \in (2, 2 + \frac{4s}{N})$, and combined powers $g(t)=|t|^{q-2} t \pm |t|^{r-2} t$, with $2\leq r < q < 2 + \frac{4s}{N}$. 
Other physical models can be found for example in the saturation effect in nonlinear optics for photorefractive media, e.g. \label{pag_comm_physics2}
$$g(t)= \frac{t^3}{1+ t^2}, \quad G(t)= \frac{1}{2} \left(t^2 - \log(1+t^2)\right)$$
(see \cite{DLQ,WCCYZH,MMP,HLRZ}).

The solutions to \eqref{problem} can be characterized as critical points of the $C^1$ functional $\mc{L} : H^s_r(\R^N) \to \R$
$$ \mathcal{L}(u) =\half \int_{\R^N} |(-\Delta)^{s/2} u|^2 - \int_{\R^N} G(u) $$
constrained on the sphere 
$$ \mathcal{S}_m = \big\{ u \in H^s_r(\R^N) \mid \|u\|^2_2= m \big\}, $$
where $\norm u_q = \left(\intRN \abs{u}^q\, dx\right)^{1/q}$ for any $q\in [1,\infty)$. 
We recall that 
$$ H^s_r(\R^N) = \big\{ u \in H^s(\R^N) \mid u(x)= u(|x|) \big\} $$
is the subspace consisting of radial functions of the Sobolev space
$$ H^s(\R^N) =\left \{u \in L^2(\R^N) \, \middle | \, \int_{\R^N} |\xi|^{2s} |\mc{F} (u)|^2 d\xi < \infty \right \}$$ 
endowed with the norm 
\begin{equation}\label{norma}
\|u\|^2_{H^s(\R^N)} = \int_{\R^N} \big(|\xi|^{2s} |\mc{F}(u)|^2 + |\mc{F} (u)|^2\big) dx
\end{equation}
where $\mc{F}$ denotes the Fourier transform.
We also recall that 
$$ (-\Delta)^{s} u = \mc{F}^{-1} \big(|\xi|^{2s}\mc{F}(u)\big) $$
for $|\xi|^{s}\mc{F}(u) \in L^2(\R^N)$ and thus $(-\Delta)^{s} u$ is a real function for a real valued $u \in L^2(\R^N)$ (see \cite{DPV}).
By Plancherel's Theorem, we have
$$ \int_{\R^N} |(-\Delta)^{s/2} u|^2 dx =
\int_{\R^N} |\mc{F}((-\Delta)^{s/2} u)|^2 d \xi =
\int_{\R^N} |\xi|^{2s} |\mc{F} (u)|^2d \xi $$
for any $u \in H^s(\R^N)$.
Therefore
$$ \|u\|^2_{H^s(\R^N)}= \int_{\R^N} \big(|(-\Delta)^{s/2} u|^2 + u^2\big) dx.$$
In \cite{Lio0}, Lions proved the following result in a radial fractional setting.
\begin{lemma}\label{compact}
Let $N \geq 2$. The space $H^s_r(\R^N)$ is compactly embedded into $L^{q+1}(\R^N)$ for all $q \in (1, 1+ \frac{4s}{N-2s})$.
\end{lemma}
However, as shown in \cite{CO} (see also \cite{BGMMV}), a result in the spirit of Radial Lemma by Strauss is not available in a fractional framework for general $0 <s \leq \frac{1}{2}$. 

\smallskip
In the present work we consider a Lagrange formulation of the problem \eqref{problem} as in \cite{HT} (see also \cite{CT}). 
For technical reasons we write $\mu = e^\lambda$ with $\lambda \in \R$ and define the $C^1$ functional $\mc{I} : \R \times H^s_r(\R^N) \to \R$ by setting 
\begin{equation}\label{functlag}
\mc{I}(\lambda, u)=\half \int_{\R^N} |(-\Delta)^{s/2} u|^2 - \int_{\R^N} G(u)\, + 
\frac{e^\lambda}{2} \bigl( \|u\|_2^2 -m \bigr).
\end{equation}
We seek for critical points $(\lambda,u) \in \R \times H^s_r(\R^N)$ of 
$\mc{I}$, namely weak solutions of $\partial_u \mc{I}(\lambda, u)=0$ and $\partial_\lambda \mc{I}(\lambda, u)=0$ or equivalently 
$$
\quad \parag{
\int_{\R^N}(-\Delta)^{s/2} u \ (-\Delta)^{s/2} \phi + e^\lambda u \phi &=& \int_{\R^N} g(u) \phi, \quad \forall \phi \in H_r^s(\R^N), \cr 
\int_{\R^N} u^2 dx &=& m.
}
$$
We implement a mini-max approach to detect normalized solutions in the nonlocal framework using a Pohozaev type function.
More precisely, inspired by the Pohozaev identity, for any $s \in (0,1)$ we introduce the Pohozaev function 
$\mc{P}:\R \times H^s_r(\R^N) \to \R$ by setting 
$$
\mc{P}(\lambda, u)=
\frac{N-2s}{2} \int_{\R^N} |(-\Delta)^{s/2} u|^2 + N \int_{\R^N} \left( \frac{e^\lambda}{2} u^2 - G(u)\right) 
$$
and the set
$$ \Omega =\big\{(\lambda,u) \in \R \times H^s_r(\R^N) \mid \mc{P}(\lambda,u)>0\big\} \cup\big\{(\lambda,0) \mid \lambda \in \R \big\}. $$
We note that, for each $\lambda\in\R$, the set
$\{u\in H_r^s(\R^N)\,|\, \mc{P}(\lambda,u)>0\}\cup\{ 0\}$ is a neighbourhood of $u=0$, and thus 
$$ \partial \Omega =\big\{(\lambda,u) \in \R \times H^s_r(\R^N) \mid \mc{P}(\lambda,u)=0, \ u \neq 0 \big\}. $$
Therefore $(\lambda, u) \in \partial \Omega$ if and only if $u \neq 0$ and $u$ satisfies the Pohozaev identity. However we emphasize that under assumptions (g1)--(g3), 
if $u \in H^s(\R^N)$ solves $\partial_u \mc{I}(\lambda, \cdot)=0$ with $\lambda \in \R$ fixed, then $\mc{P}(\lambda, u)=0$ when $s \in (1/2,1)$. A similar result for $s \in (0,1/2]$ is not available since the weak solutions are not $C^1$, in general (see \cite{BKS}).

In spite of this lack of regularity, which is a special feature of the nonlocal framework, we recognize a Mountain Pass structure \cite{AR} for the functional $\mc{I}$, where the mountain is given by the subset $\partial\Omega$. We refer to it as the \emph{Pohozaev's mountain}. 
This approach can be useful to deal with different problems in other contexts.

Inspired by \cite{HT,IT}, we need to use a new variant of the Palais-Smale condition which takes into account the Pohozaev identity, and we establish some deformation theorems which enable us to perform our minimax arguments in the product space $\R \times H^s_r(\R^N)$. 

As a byproduct, our solutions satisfy the Pohozaev identity, even if we assume that $f$ is a continuous function (see Corollary \ref{coroll_esist_Pm}).
We also note that solutions with the Pohozaev identity are essential in the following sense: our deformation argument shows that critical points with the Pohozaev identity just contribute to the topology. 
That is, solutions without the Pohozaev identity are deformable with a suitable deformation flow and have no topological contribution.

Firstly we prove the following existence results for \eqref{problem}. 

\begin{theorem}\label{S:1.1}
Suppose $N\geq 2$ and \textnormal{(g1)--(g3)}.
Then there exists $m_0 \geq 0$ such that for any $m>m_0$, the problem \eqref{problem} has a solution.
\end{theorem}	
	
\begin{theorem}\label{S:1.12}
		Suppose $N\geq 2$, \textnormal{(g1)--(g3)} and 
		\begin{itemize}
		\item[\textnormal{(g4)}] $\lim_{t \to 0} \frac{g(t)}{|t|^{\frac{4s}{N}} t} = + \infty.$
		\end{itemize}
		Then for any $m >0$, the problem \eqref{problem} has a solution.
\end{theorem}
	
We highlight that the found solution is actually a minimum for $\mc{L}$ constrained to the sphere (see Proposition \ref{minimizing}), which furnishes a strong indication to its stability properties.
The techniques employed in \cite{Shi0} for the local case $s=1$, to get directly the existence of a minimum for $\mc{L}$, are not easily adaptable to the fractional framework, because of the need of a control on the tails in the Brezis-Lieb lemma and in the Concentration-Compactness techniques. 
Anyway, our method not only gets around these difficulties, but moreover it is also suitable to get multiple solutions.

\label{pag_comm_minimum}

Indeed, if we also suppose the oddness of $g$, namely
\begin{itemize}
	\item[(g5)] $g(-t)= - g (t)$ for all $t \in \R$,
\end{itemize}
we have $\mc{I}(\lambda,-u)= \mc{I}(\lambda,u)$ for all $(\lambda,u) \in \R \times H^s_r(\R^N)$ and we can establish the existence of infinitely many $L^2$ constrained standing waves solutions for the (fNLS) equation. 

We prove the following multiplicity result.

\begin{theorem}\label{S:1.13}	
Suppose $N\geq 2$ and \textnormal{(g1)--(g3)} and \textnormal{(g5)}. 
Then we have:
\begin{itemize}
	\item[\textnormal{(i)}]
	 For any $k \in \N$ there exists $m_k \geq 0$ such that for each $m > m_k$, the problem \eqref{problem} has at least $k$ nontrivial, distinct pairs of solutions. 
	\item[\textnormal{(ii)}]
	In addition assume \textnormal{(g4)}. For any $m >0$ the problem \eqref{problem} has countably many solutions $(u_n)_{n \in \N}$, which satisfy 
	$$\mathcal{L}(u_n) <0 \quad \hbox{for all $ n \in \N$},$$ 
	$$\mathcal{L}(u_n) \to 0 \quad \hbox{as $ n \to + \infty$}.$$
\end{itemize}
\end{theorem}

We remark that our subcritical multiplicity result seems new even in the case of the pure power $g(t)=|t|^{q-2}t$ and in the non-monotone case of competing powers $g(t)= |t|^{q-2}t - |t|^{r-2}t$ and it has a physical relevance since it describes the existence of multiple bound states with arbitrary high energies (see e.g. \cite{DZ}).
We stress that the analytical solutions for fractional differential equations are still limited, while there is a large amount of numerical methods in discretizing the fractional differential operators. 
In Theorem \ref{S:1.13} we furnish an analytical rigorous approach to detect infinitely many symmetric solitons, which can be applied to the computation of ground and excited states to (fNLS) equations arising from Bose-Einstein condensation theory or nonlinear optics phenomena with saturation.

\label{pag_comm_physics3}

\medskip
The paper is organized as follows. 
In Section \ref{section:2}, we establish some preliminaries related to the unconstrained problem. 
In Section \ref{section:3} we give the Lagrange formulation of the problem \eqref{problem} and a description of the geometry of an auxiliary functional in a product space. 
Section \ref{section:4} concerns with the Palais-Smale-Pohozaev ((PSP) for short) condition and Section \ref{section:5} is devoted to the construction of the deformation argument under the (PSP) condition. 
Section \ref{section:6} deals with our minimax procedure to detect the normalized solutions by means of the Pohozaev's mountain. 
Finally in Section \ref{section:7} we derive the multiplicity result of infinitely many normalized solutions when $g$ is odd.


\setcounter{equation}{0} 
\section{\label{section:2}The unconstrained problem in $H^s(\R^N)$}

In this section we consider the unconstrained fractional equation
\begin{equation} \label{problemxx}
\parag{
&(-\Delta)^{s} u + \mu u =	f(u)& \quad \hbox{in $\R^N$,} \cr 
& u \in H^s(\R^N),& \cr
& u>0,&
}
\end{equation}
where $s \in (0,1)$, $N \geq 2$, $\mu >0$ is fixed and $f$ satisfies the following assumptions
\vskip2pt
\begin{itemize}
	\item[(f1)] $f : \R \to \R$ continuous, $f(t)= 0$ for $t \leq 0$ and $\lim_{t \to 0} \frac{f(t)}{t}=0;$
	\item[(f2)] $\limsup_{|t| \to \infty} \frac{f(t)}{|t|^q} =0$ where $ q \in (1, 1+ \frac{4s}{N-2s})$;
	\item[(f3)] there exists $t_0 \in \R$, $t_0 > 0$ such that $F(t_0) >\frac{\mu}{2}t_0^2$, where $F(t)= \int_0^t f(\tau) d\tau$. 
\end{itemize}
\vskip 2pt

Under the assumptions (f1)-(f2), it is standard to show that any weak solution of \eqref{problemxx} is a critical point of the $C^1$ functional $J: H^s(\R^N) \to \R$ defined by
$$
J (u)=\half \int_{\R^N} | (-\Delta)^{s/2} u|^2 \, dx + \frac{\mu}{2} \int_{\R^N} u^2 \, dx - \int_{\R^N} F(u) \, dx. 
$$

In the celebrated paper \cite{BL1}, for the local case $s=1$, Berestycki and Lions proved the existence of a classical solution to \eqref{problemxx}, which is radially symmetric and has an exponentially decay, under the assumption (f1)--(f3). 
These conditions are almost optimal for the existence of \eqref{problemxx}. 
The found solution is of least energy among all nontrivial solutions. 
Indeed, when $s=1$, in \cite{JT1} Jeanjean and the third author proved that the last energy solution is indeed a Mountain Pass (MP for short) solution. 
Successively Byeon, Jeanjean, Maris \cite{BJM} showed that every least energy solution of \eqref{problemxx} is radially symmetric, up to translations.

For the nonlocal case $s \in (0,1)$, the equivalence between MP weak solution and least energy solution is still a partially open problem.

We begin to recall that in the recent paper \cite{BKS}, Byeon, Kwon and Seok established the following results.
\begin{proposition}\label{Regularity}
	Suppose \textnormal{(f1)-(f2)}. Let $u \in H^s(\R^N)$ be a weak solution 
	of the fractional equation \eqref{problemxx}.
	Then $u \in C^1(\R^N)$ if one of the following assumption holds:
	\begin{itemize}
		\item[\textnormal{(i)}]	$s \in (1/2,1)$;
	\item[\textnormal{(ii)}]
	$s \in (0,1/2]$ and $f \in C^{0,\alpha}_{loc}(\R)$ for some $\alpha \in (1-2s,1)$.
	\end{itemize}	
\end{proposition}

\begin{proposition}\label{Pohozaev-prop}
	Suppose \textnormal{(f1)--(f3)} and
	\begin{itemize}
		\item[\textnormal{(f4)}] if $s \in (0, 1/2]$, $f \in C^{0, \alpha}_{loc} (\R)$ for some $\alpha \in (1-2s,1)$.
	\end{itemize}
	\vskip 2pt
	Then every weak solution $u \in H^s(\R^N)$ of the fractional equation of \eqref{problemxx} satisfies the Pohozaev identity
	\begin{equation}\label{Pohozaev}
	\frac{N-2s}{2} \int_{\R^N} | (-\Delta)^{s/2}u|^2 \, dx + N \int_{\R^N} \left( \frac{\mu}{2} u^2 - F(u)\right) dx =0.
	\end{equation}
\end{proposition}
We remark that the $C^1$ regularity of the weak solution seems crucial for proving a Pohozaev type identity.
Under (f1)--(f3) we know \cite{BKS} that each weak solution of \eqref{problemxx} belongs to $H^s(\R^N) \cap C^\beta (\R^N)$ with $\beta \in (0, 2s)$ and thus it is not known if the Pohozaev identity holds when $s \in (0,1/2]$, without additional regularity assumptions on the nonlinearity $f$.

 In \cite[Theorem 1.2]{BKS}, they also investigated the existence of MP weak solutions of \eqref{problemxx}. 
 We recall that a weak solution $u$ is said of MP type if 
\begin{equation}\label{MPlevel}
J(u) = C_{mp},
\end{equation}
where
$$
C_{mp} = \inf_{\gamma\in\Gamma}\max_{t\in [0,1]} J(\gamma(t))
$$
and
\begin{equation}\label{classpath}
\Gamma = \big\{ \gamma(t)\in C\big([0,1], H^s(\R^N)\big) \mid \gamma(0)=0, \, J(\gamma(1))<0\big\}.
\end{equation}
As for $s=1$, the functional $J$ does not satisfies the Palais-Smale condition at level $C_{mp}$ under the assumptions (f1)--(f3), thus one can not directly apply the MP theorem. For the local case $s=1$, any weak solution is $C^1$ and it satisfies the Pohozaev identity, so that one can reduce the search of MP solutions to that of minimizers on the Pohozaev type constraint.

For the fractional case, this approach seems to work for $s \in (1/2, 1)$, while requires additional regularity on the nonlinearity if $s \in (0,1/2]$.
 
Conversely in \cite{BKS}, the authors established that every minimizer of $J$ on the Pohozaev type constraint corresponds to a MP weak solution and derived some radially symmetric properties of the minimizer using a fractional version of the Polya-Szego inequality.
Namely they introduce the Pohozaev functional $\mathcal{P}: H^s(\R^N) \to \R$ by setting 
$$ \mathcal{P}(u) = \frac{N-2s}{2} \int_{\R^N} |(-\Delta)^{s/2}u|^2 + N\int_{\R^N} \left(\frac{\mu}{2} u^2 - F(u)\right) $$
and 
$${P}=\big\{ u \in H^s (\R^N) \setminus \{ 0 \} \, | \, \mathcal{P}(u)=0 \big \},$$
$$C_{po} = \min_{u \in{P}} I(u).$$
In \cite[Theorem 1.2]{BKS} they established the following result.
\begin{theorem}\label{Byeonminimizers}
Assume \textnormal{(f1)--(f3)}. Fix $s \in (0,1)$ and $\mu >0$. Then 
	\begin{itemize}
		\item[\textnormal{(i)}]
			there exists a minimizer of $J$ subject to ${P}$;
			\item[\textnormal{(ii)}]
			every minimizer of $J$ subject to ${P}$ is a MP weak solution of 
		\eqref{problemxx};
		\item[\textnormal{(iii)}]	
			every minimizer of $J$ subject to ${P}$ is radially symmetric up to a translation.
\end{itemize}
\end{theorem}
From Theorem \ref{Byeonminimizers} it follows that 	
$$ C_{mp} =C_{po}. $$
However the equivalence between Mountain Pass solutions and least energy solutions is shown for $s\in(1/2,1)$, while it is yet an open problem for $s\in(0,1/2]$ under the assumptions (f1)--(f3). 
In \cite{BKS}, this equivalence is established under the same regularity assumption of Proposition 2.1, namely $f\in C^{0,\alpha}(\R^N)$ for some $\alpha\in (1-2s,1)$.
In the following sections, in contrast, we will show that under $L^2$ constraint least energy solutions have Mountain Pass characterization.
See Proposition \ref{minimizing}. \label{page_minim_2}


\setcounter{equation}{0} 
 \section{Lagrange formulation and geometry of $\mc{I}(\lambda,u)$} \label{section:3}

We come back to the constrained case and we consider the Lagrange formulation of the problem \eqref{problem} in the space of radially symmetric functions $H_r^s(\R^N)$.
Namely, we seek for critical points of the functional $\mc{I}: \R \times H^s_r(\R^N) \to \R$ defined in \eqref{functlag}, i.e.
\begin{equation}\label{equation}
 \mc{I}(\lambda,u)=\half \int_{\R^N} |(-\Delta)^{s/2} u|^2 - \int_{\R^N} G(u)\, + 
 \frac{e^\lambda}{2} (\|u\|_2^2 -m).
 \end{equation}
Under the assumption (g1)--(g3), it is standard to prove that $\mc{I}$ is $C^1$ in the product space $\R \times H^s_r(\R^N)$.

It is immediate to recognize that for any $m >0$
$$\mc{I}(\lambda,u)=\mc{J}(\lambda,u) -{e^\lambda\over 2}m$$
where $\mc{J}:\R \times H^s_r(\R^N) \to \R$ is the $C^1$ functional defined by 
$$
\mc{J}(\lambda, u)=\half \int_{\R^N} |(-\Delta)^{s/2} u|^2 - \int_{\R^N} G(u)\, + 
\frac{e^\lambda}{2} \int_{\R^N} u^2.
$$
For a fixed $\lambda \in \R$, $u$ is critical point of $\mc{J}(\lambda, \cdot)$ means that $u$ solves, in the weak sense,
\begin{equation}\label{problemxxx}
\parag{
&(-\Delta)^{s} u + e^\lambda u = g(u)& \qquad \hbox{in $\R^N$,} \\ 
& u \in H^s_r(\R^N).&
}
\end{equation}
Inspired by the Pohozaev identity, for any $s \in (0,1)$ we also introduce the Pohozaev functional $\mc{P}:\R \times H^s_r(\R^N) \to \R$ by setting 
$$
\mc{P}(\lambda, u)=
\frac{N-2s}{2} \int_{\R^N} |(-\Delta)^{s/2} u|^2 + N \int_{\R^N} \left( \frac{e^\lambda}{2} u^2 - G(u)\right).
$$
By Proposition \ref{Pohozaev-prop}, it follows that for any $\lambda \in \R$, if $u \in H^s_r(\R^N)$ solves \eqref{problemxxx}, then $\mc{P}(\lambda, u)=0$ when $s \in (1/2,1)$.
A similar result for $s \in (0,1/2]$ is not known under (g1)--(g3).

\smallskip

Now set
$$ \Omega =\big\{(\lambda,u) \in \R \times H^s_r(\R^N) \mid \mc{P}(\lambda,u)>0\big\} \cup\big\{(\lambda,0) \mid \lambda \in \R \big\}.$$
We note that for each $\lambda\in\R$, $\mc{P}(\lambda,u)>0$ in a small neighbourhood of $u=0$ except $0$. Thus we have 
$$ \partial \Omega =\big\{(\lambda,u) \in \R \times H^s_r(\R^N) \mid \mc{P}(\lambda,u)=0, \ u \neq 0 \big\}, $$
which we call the \emph{Pohozaev mountain} for $\mc{J}(\lambda,u)$.
We remark that $(\lambda, u) \in \partial \Omega$ if and only if $u \neq 0$ and $u$ satisfies the Pohozaev identity $\mc{P}(\lambda, u)=0$.

Set
\begin{equation}\label{muzero}
\mu_0 = 2 \sup_{s \in \R, s \neq 0} {G(s)\over s^2}, 
\end{equation}
we can deduce $\mu_0 \in (0, \infty]$ under the assumptions (g1)--(g3). 
In what follows, we denote 
\begin{equation}\label{lambdazero}
\lambda_0 = \log \mu_0, \quad \hbox{if} \ \mu_0 \in (0, \infty),
\end{equation}
otherwise $\lambda_0 = + \infty$.

Taking into account that $$1+ \frac{4s}{N} < 1 + \frac{4s}{N-2s},$$ we deduce by Theorem \ref{Byeonminimizers} that for any $\lambda \in (- \infty, \lambda_0)$ the functional 
$$ u \mapsto \mc{J}(\lambda,u);\, H^s_r(\R^N) \to \R $$
has a minimizer $u_\lambda$ subject to 
$$(\partial \Omega)_{\lambda} =
\big\{ u \in H_r^s (\R^N) \setminus \{ 0 \} \mid \mc{P}(\lambda,u)=0 \big\}, $$
namely
\begin{equation}\label{minimumpohozaev}
\mc{J}(\lambda,u_\lambda)= \min_{u \in (\partial \Omega)_{\lambda}} \mc{J}(\lambda,u).
\end{equation}
Furthermore by $(ii)$ of Theorem \ref{Byeonminimizers} such $u_\lambda$ is a Mountain Pass critical point of $\mathcal{J}(\lambda, \cdot)$ at level $a(\lambda)$, i.e.
$$ \mc{J}(\lambda,u_\lambda)= a(\lambda) $$
where 
\begin{equation}\label{MPlevelx}
a(\lambda) = \inf_{\gamma\in\Gamma(\lambda)}\max_{t\in [0,1]} \mc{J}(\lambda,\gamma(t))
 \end{equation}
and
\begin{equation}\label{classpathx}
\Gamma(\lambda) = \big\{ \gamma \in C\big([0,1], H^s_r(\R^N)\big) \mid \gamma(0)=0, \, \mc{J}(\lambda,\gamma(1))<0\big\}.
\end{equation}
We notice that $\lambda \mapsto a(\lambda);\, (-\infty,\lambda_0)\to\R$ is strictly monotone increasing on $\R$.

\begin{lemma}\label{lem_buona_def}
Let $\lambda \in \R$. 
Then the following statements are equivalent:
\begin{itemize}
\item[\textnormal{(a)}] $\lambda < \lambda_0$.
\item[\textnormal{(b)}]
There exists a $t_0=t_0(\lambda)>0$ such that
$$ G(t_0) > \frac{e^{\lambda}}{2} |t_0|^2.$$
\item[\textnormal{(c)}]
There exists $u\in H^s_r(\R^N)\setminus \{0\}$ such that $\mc{P}(\lambda, u)=0$; in particular $(\partial \Omega)_{\lambda}\neq \emptyset$.
\item[\textnormal{(d)}]
$\Gamma(\lambda)\neq \emptyset$, and thus $a(\lambda)$ is well defined.
\end{itemize}
As further consequence, we see that $\partial \Omega \neq \emptyset$. Finally, $a(\lambda)>0$.
\end{lemma}

\claim Proof.
(a) $\iff$ (b).
This is a straightforward consequence of the definition of $\lambda_0$.

(b) $\implies$ (c)
Let $u\in H^s_r(\R^N)$ to be fixed. 
We have, for $t>0$,
$$\mc{P}(\lambda, u(\cdot/t))= \frac{N-2s}{2} t^{N-2s} \int_{\R^N} |(-\Delta)^{s/2}u|^2 - N t^{N} \int_{\R^N} \left(G(u) -\frac{ e^{\lambda}}{2} u^2 \right).$$
We notice that $\mc{P}(\lambda, u(\cdot/t))>0$ for small $t>0$. In order to get a $\bar{t}$ such that $\mc{P}(\lambda, u(\cdot/\bar{t}))=0$ we need the quantity
$$ \int_{\R^N} \left(G(u) - \frac{e^{\lambda}}{2} u^2\right)$$
to be positive. 
For any $R>0$ we choose a smooth $u=u_R\in C^{\infty}_c$ such that $u_R=t_0$ in $B_R(0)$ and $u_R=0$ out of $B_{R+\frac{1}{R^N}}(0)$, $0\leq u_R\leq t_0$. 
We set
$$C=\sup_{t\in[0,t_0]} \left|G(t)-\tfrac{e^{\lambda}}{2} |t|^2\right|<+\infty.$$
Then
\begin{eqnarray*}
\lefteqn{ \int_{\R^N} \left(G(u_R) - \frac{e^{\lambda}}{2} u_R^2 \right)} \\
&=& \int_{B_{R+\frac{1}{R^N}}\setminus B_R}\left(G(u_R) - \frac{e^{\lambda}}{2} u_R^2 \right) +\int_{B_R} \left(G(u_R) - \frac{e^{\lambda}}{2} u_R^2 \right) \\
&\geq& -C \abs{B_{R+\frac{1}{R^N}} \setminus B_R} + |B_R| \left(G(t_0) - \frac{e^{\lambda}}{2} |t_0|^2\right) \to +\infty
\end{eqnarray*}
and in particular it is positive for a sufficiently large $R$.

(c) $\implies$ (d).
Let $u\in H^s_r(\R^N) $, $u\nequiv 0$ such that $\mc{P}(\lambda,u)=0$. 
We define $\gamma(t)=u(\cdot/t)$ for $t\neq 0$ and $\gamma(0)= 0$, so that $\gamma:[0, \infty)\to H^s_r(\R^N)$ is continuous.
We have
$$\mc{J}(\lambda, \gamma(t))= \frac{1}{2} t^{N-2s} \int_{\R^N} |(-\Delta)^{s/2}u|^2 - t^N \int_{\R^N} \left(G(u) - \frac{e^{\lambda}}{2} u^2\right).$$
Noting $\int_{\R^N}\left(G(u)-\frac{e^\lambda}2 u^2\right)>0$ by $\mc{P}(\lambda,u)=0$, we have
$\mc{J}(\lambda,\gamma(t))\to -\infty$ as $t\to\infty$ and thus $\Gamma(\lambda)\not=\emptyset$.

(d) $\implies$ (b).
If $\gamma\in \Gamma(\lambda)$, then $\mathcal{J}(\lambda,\gamma(1))<0$, thus
$$\int_{\R^N} \left(G(\gamma(1)) - \frac{e^{\lambda}}{2} \gamma(1)^2\right)>0,$$
which implies that there exists an $x_0\in \R^N$ such that
$$G(\gamma(1)(x_0)) - \frac{e^{\lambda}}{2} \gamma(1)^2(x_0)>0.$$
The claim comes setting $t_0=\gamma(1)(x_0)$.

Finally, by Theorem \ref{Byeonminimizers}, there exists a Pohozaev minimum $u_{\lambda}$ which is also a Mountain Pass solution, thus $\mathcal{J}(\lambda,u_{\lambda})=a(\lambda)$, $D_u \mathcal{J}(\lambda, u_{\lambda})=0$ and $\mathcal{P}(\lambda,u_{\lambda})=0$, which imply
$$ a(\lambda)=\frac{s}{N} \| (-\Delta)^{s/2} u_\lambda \|^2_2>0. 
\QED
$$


\begin{remark} 
Assume $\lambda_0 <+\infty$. 
We observe that, in this case, for $\lambda\geq \lambda_0$ we have $\mc{P}(\lambda, u)\geq 0$ and $\mc{J}(\lambda, u)\geq 0$ for each $u$, both strictly positive for $u\nequiv 0$. 
This means that $[\lambda_0,+\infty)\times H_s^r(\R^N) \subset \Omega$.
\end{remark}

\begin{proposition}\label{S:2.33}
	Assume \textnormal{(g1)--(g3)} and $\lambda_0 \in \R$. We have 
	\smallskip
	\begin{itemize}
		\item[\textnormal{(a)}] if $(\lambda,u) \in \partial \Omega$ for some $u \in H^s_r(\R^N)$, then $\lambda < \lambda_0$.
		\item[\textnormal{(b)}] $\lim_{\lambda \to \lambda_0^-} a(\lambda)=+\infty$.
	\end{itemize}
\end{proposition}

\claim Proof.
Let $(\lambda, u) \in \partial \Omega$, namely $\mathcal{P}(\lambda,u)=0$ and $u \neq 0$. 
This implies that for some $x \in\R^N$ 
$$ G(u(x)) - \frac{e^\lambda}{2} u(x)^2 >0 $$
and thus $\lambda < \lambda_0$ and (a) holds.

Now we show point (b). Let $\lambda<\lambda_0$; by contradiction, since by Lemma \ref{lem_buona_def} $a(\lambda)$ is increasing and strictly positive, we assume that $a(\lambda)\to c \in (0, +\infty)$ as $\lambda \to \lambda_0^-$, from which we deduce that $\norm{(-\Delta)^{s/2}u_{\lambda}}_2$ is bounded. Moreover, for any fixed $\delta>0$ there exists a suitable $C_{\delta}>0$ such that
$$G(s)\leq \delta |s|^2 + C_{\delta} |s|^{p+1},$$
where we recall that $p=1+{4s\over N}$.
 
Thus we have by the fractional Gagliardo-Nirenberg inequality \cite{Par0} (see also \cite{BGMMV}) and the fact that $\norm{(-\Delta)^{s/2} u_{\lambda}}_2$ is bounded,
\begin{eqnarray*}
	0&=& \frac{1}{N} \mathcal{P}(\lambda,u_{\lambda}) \geq 
	\frac{N-2s}{2N}
	\norm{(-\Delta)^{s/2}u_{\lambda}}_2^2 +
	\left(\frac{e^{\lambda}}{2} -\delta\right) \norm{u_{\lambda}}_2^2 - C_{\delta}\norm{u_{\lambda}}_{p+1}^{p+1} \\
	&\geq& \frac{N-2s}{2N} \norm{(-\Delta)^{s/2}u_{\lambda}}_2^2 + \left(\frac{e^{\lambda}}{2} -\delta\right) \norm{u_{\lambda}}_2^2 
 - C'C_{\delta}
	 \norm{(-\Delta)^{s/2} u_{\lambda}}_2^2 \norm{u_{\lambda}}_2^{p-1} \\
	&\geq & \left(\frac{e^{\lambda}}{2} -\delta\right) \norm{u_{\lambda}}_2^2 - C'' C_{\delta} \norm{u_{\lambda}}_2^{\frac{4s}{N}}
\end{eqnarray*}
for some $C'$, $C''>0$. 
By choosing $\delta < \frac{e^{\lambda}}{2}$, since $\frac{4s}{N}<2$, also $\norm{u_{\lambda}}_2$ must be bounded, which means that $(u_\lambda)_{\lambda<\lambda_0}$ is bounded in $H^s_r(\R^N)$. 
Hence, up to a subsequence, $u_{\lambda}\rightharpoonup u_0$ in $H^s_r(\R^N)$. By Lemma \ref{compact} and taking into account that $\partial_u\mathcal{J}(\lambda, u_{\lambda})=0$, we deduce that $u_\lambda \to u_0$ strongly in $H^s_r(\R^N)$ with $\mathcal{J}(\lambda_0, u_0)=c$, $\partial_u\mathcal{J}(\lambda_0, u_0)=0$, $\mathcal{P}(\lambda_0, u_0)=0$.
Since $c>0$, we have $u_0 \neq 0$. By $\mathcal{P}(\lambda_0, u_0) =0$, we conclude
$$ G(u_0(x)) - \frac{e^{\lambda_0}}{2} u_0(x)^2 >0 $$
for some $x \in \R^N$, which contradicts the definition of $\lambda_0$. 
\QED

\bigskip

In the next result, we consider the case $\lambda_0= + \infty$ and we investigate the behaviour of $a(\lambda)$ for $\lambda$ large.
\begin{proposition}\label{S:lim}
Assume that $\lambda_0= + \infty$. Then
$$ \lim_{\lambda\to +\infty} {a(\lambda)\over e^\lambda} = +\infty.$$
\end{proposition}

\claim Proof.
By (g1)-(g2) we have that for any $\delta>0$ there exists $C_\delta>0$ such that for all $s\in\R$
\begin{eqnarray}
\abs{g(s)} &\leq& \delta\abs s^p + C_\delta\abs s, \label{maggiorazione-0} \\
\abs{G(s)} &\leq& 
\frac{\delta}{p+1}|s|^{p+1} + \frac{C_\delta}{2} \abs s^2, \label{maggiorazione} 
\end{eqnarray}
where $p = 1 + 4s/N$. 
We also denote by $b(\delta)$ the MP value of $\mc{H}_{\delta} : H_r^s(\R^N) \to \R$ defined by
$$\mc{H}_{\delta}(v) = \half \norm{(-\Delta)^{s/2}v}_2^2 + \half \norm{v}_2^2 -\frac{\delta}{p+1}\norm{v}_{p+1}^{p+1}.$$
It is easy to see that
$$b(\delta) \to +\infty \quad \hbox{as}\ \delta\to 0^+.$$
For $v\in H^s_r(\R^N)\setminus\{ 0\}$, we set 
$$u_\theta(x)=\theta^{N/2}v(\theta x),$$
and for simplicity we write $\mu = e^\lambda$ and $\mc{J}(\mu,\cdot) =\mc{J}(\lambda,\cdot)$. 
By \eqref{maggiorazione}, we pass to evaluate
$$\mc{J}(\mu,u_\theta) \geq \theta^{2s} \left(\half \norm{(-\Delta)^{s/2}v}_2^2 
 + \half (\mu - C_\delta) \theta^{-2s} \norm{v}_2^2 
 - \frac{\delta}{p+1} \norm{v}_{p+1}^{p+1} \right).$$
Setting $\theta = (\mu - C_\delta)^{1/{2s}}$ for $\mu > C_\delta$, we have
\begin{equation}\label{key1}
\mc{J}(\mu, u_{(\mu - C_\delta)^{1/{2s}}} ) \geq (\mu -C_\delta) \mc{H}_{\delta} (v)
\end{equation}
and hence
\begin{equation}\label{key2}
{\mc{J}(\mu, u_{(\mu - C_\delta)^{1/{2s}}})\over\mu} \geq \frac{\mu - C_\delta}{\mu}
\mc{H}_\delta (v).
\end{equation}
Thus we have
\begin{equation}\label{key3}
 {a(\mu)\over \mu} \geq {\mu-C_\delta\over \mu} b(\delta);
\end{equation}
since $\delta>0$ is arbitrary, we derive
$$ \lim_{\mu\to +\infty} {a(\mu)\over\mu} = +\infty. 
\QED
$$

\medskip

\begin{proposition}\label{out}
	Assume \textnormal{(g4)} in addition to \textnormal{(g1)--(g3)}. 
	Then
	\begin{equation}\label{third}
	\lim_{\lambda\to -\infty} {a(\lambda)\over e^{\lambda}} = 0.
	\end{equation}
\end{proposition}
\claim Proof.
We fix $u\in H^s_r(\R^N) \cap L^\infty(\R^N)$ with $\norm u_\infty=1$. Set $p = 1+ \frac{4s}{N}$, there exists $M_r>0$
such that for all $r \in (0,1]$ 
$$G(ru(x)) \geq \frac{M_r}{p+1}r^{p+1} \abs{u(x)}^{p+1}, \ \quad \forall x\in\R^N$$
with 
$$M_r\to+\infty \quad \hbox{as}\ r\to 0.$$
We write again $\mu=e^{\lambda}$ for the sake of simplicity. Therefore for $t>0$ we have 
\begin{eqnarray*}
\lefteqn{ \mc{J}(\mu,r u(x/t)) \leq \half r^2 t^{N-2s} \|(-\Delta)^{s/2}u \|_2^2 +{\mu\over 2}r^2t^N \norm u_2^2
	-\frac{M_r}{p+1} r^{p+1} t^{N} \|u\|_{p+1}^{p+1} } \\ 
&=& r^2 \mu^{-\frac{N-2s}{2s}} \left( \half t^{N-2s} \mu^{\frac{N-2s}{2s}} 
	\|(-\Delta)^{s/2}u \|_2^2 +{1\over 2} \mu^{\frac{N}{2s}} t^N \norm u_2^2
	-\frac{M_r}{p+1} r^{\frac{4s}{N}} \mu^{\frac{N-2s}{2s}} t^{N} \|u\|_{p+1}^{p+1}\right) 	\\ 
&=& r^2 \mu^{-\frac{N-2s}{2s}} \left( {1\over 2}\tau^{N-2s} 
	\|(-\Delta)^{s/2}u \|_2^2 +{1\over 2} \tau^{N} \norm u_2^2
	-\frac{M_r}{p+1} r^{\frac{4s}{N}} \mu^{-1} \tau^{N} \|u\|_{p+1}^{p+1} \right)
	\end{eqnarray*}
after setting $\tau = \mu^{\frac{1}{2s}}t$. Moreover choosing $r=\mu^{\frac{N}{4s}}$ we infer 	
\begin{eqnarray*}
\lefteqn{ \mc{J}\left(\mu,\mu^{\frac{N}{4s}} u(\cdot/(\mu^{-1/(2s)}\tau))\right) }
\\ &\leq& \mu\left( \half\tau^{N-2s} 	\|(-\Delta)^{s/2}u \|_2^2 + \half \tau^N \norm u_2^2
-\frac{M_{\mu^{N/(4s)}}}{p+1} \tau^N \|u\|_{p+1}^{p+1}\right). 
\end{eqnarray*}
For $\mu\in (0,1)$, the map
$$ \tau \mapsto \mu^{\frac{N}{4s}} u(\cdot/\mu^{-1/(2s)}\tau);\, (0,\infty)\to H^s_r(\R^N)$$
can be regarded as a path in $\Gamma(\lambda)$ after rescaling. 
Thus
$$ {a(\mu)\over \mu} \leq \max_{\tau\in [0,\infty)} 
\left( \half 
\|(-\Delta)^{s/2}u \|_2^2 \tau^{N-2s} + \half \norm u_2^2\tau^N - \frac{M_{\mu^{N/(4s)}}}{p+1} \|u\|_{p+1}^{p+1}\tau^{N}\right). $$
Since $M_{\mu^{N/(4s)}} \to \infty$ as $\mu\to 0$, we derive the conclusion. 
\QED

\medskip

\begin{proposition}\label{S:2.2}
	Assume \textnormal{(g1)--(g3)}. 
	Then we have 
\begin{itemize}
\item[\textnormal{(a)}] $\mc{J}(\lambda,u)\geq 0$ for all $(\lambda,u)\in \Omega$;
\item[\textnormal{(b)}] $\mc{J}(\lambda,u)\geq a(\lambda)>0$ for all $(\lambda,u)\in \partial\Omega$.
\end{itemize}
\end{proposition}

\claim Proof.
We notice that for all $(\lambda,u)\in \Omega$ 
$$ \mc{J}(\lambda,u) \geq \mc{J}(\lambda,u) - \frac{\mc{P}(\lambda,u)}{N} 
= \frac{s}{N} \norm{(-\Delta)^{s/2} u}_2^2 \geq 0 $$
and thus (a) follows.

The proposition (b) follows from the fact that every minimizer of $\mc{J}(\lambda,\cdot)$ subject to $(\partial \Omega)_{\lambda}$ is a mountain pass weak solution of \eqref{problemxx} at level $a(\lambda)$ (see (ii) of Theorem \ref{Byeonminimizers}).
\QED

\bigskip

We are ready to show that for any $m >0$ the functional $\mc{I}$ is bounded from below on the Pohozaev set $\partial \Omega$.
\begin{proposition}\label{negativo}
Assume \textnormal{(g1)--(g3)}.
For any $m>0$, we set
$$B_m=\inf_{\lambda < \lambda_0} \left(a(\lambda) -\frac{e^\lambda}{2}m \right)$$
and
$$B'_m= \inf_{(\lambda,u)\in\partial\Omega}\mc{I}(\lambda,u).$$
Then
\begin{equation}\label{limitinf}
B'_m \geq B_m >-\infty.
\end{equation}
\end{proposition}
\claim Proof.
Let $m>0$. 
If $(\lambda, u) \in \partial \Omega$, by (b) of Proposition \ref{S:2.2} we have
$$\mathcal{I}(\lambda, u) = \mathcal{J}(\lambda, u) - \frac{e^{\lambda}}{2} m \geq a(\lambda) - \frac{e^{\lambda}}{2} m;$$
since, by (a) of Proposition \ref{S:2.33} it results that $\lambda < \lambda_0$, we have, passing to the infimum,
$$B'_m \geq B_m.$$
We distinguish now two cases.
\smallskip
Firstly we assume $\lambda_0 \in \R$. From (b) of Proposition \ref{S:2.33}
we have $a(\lambda) \to + \infty$ as $\lambda \to \lambda_0^-$, and thus we conclude 
$$\inf_{\lambda < \lambda_0} \left(a(\lambda) -\frac{e^\lambda}{2}m \right)> -\infty.$$
Secondly, we suppose that $\lambda_0= +\infty$. We have
$$a(\lambda) -\frac{e^\lambda}{2}m = e^{\lambda} \left(\frac{a(\lambda)}{e^{\lambda}}- \frac{m}{2}\right)$$
and thus, by Proposition \ref{S:lim}
$$\inf_{\lambda \in \R} \left(a(\lambda) -\frac{e^\lambda}{2}m \right)> -\infty.
\QED
$$


\setcounter{equation}{0} 
\section{Palais-Smale-Pohozaev condition}\label{section:4}

Firstly we introduce the notations:
\begin{eqnarray*}
K_b &=& \big\{ (\lambda,u)\in \R\times H^s_r(\R^N) \mid \mc{I}(\lambda, u)=b,\,
 \partial_\lambda\mc{I}(\lambda, u)=0,\, \partial_u\mc{I}(\lambda, u)=0 \big\}, \\
K^{PSP}_b &=& \big\{ (\lambda,u)\in \R\times H^s_r(\R^N) \mid (\lambda,u) \in K_b, \,
\mc{P}(\lambda,u)=0 \big\}.
\end{eqnarray*}
Clearly, we have $K^{PSP}_b \subset K_b.$
We note that for the definition of $K^{PSP}_b$ we do not need additional regularity about $g$.

Under the assumptions (g1)--(g3), it seems difficult to verify the standard Palais-Smale condition for the functional $\mc{I}$.
Therefore we cannot recognize that the set $K_b$ is compact.

Inspired \cite{HT,IT}, we introduce the Palais-Smale-Pohozaev (shortly (PSP)) condition, which is a weaker compactness condition than the standard Palais-Smale one. Such (PSP) condition takes into account the scaling properties of $\mc{I}$ through the Pohozaev functional $\mc{P}$. 
Using this new condition we will show that $K^{PSP}_b$ is compact when $b<0$.


\subsection{(PSP) condition}

We give the definition of (PSP) condition in the radial setting.

\smallskip
\begin{definition}\label{PSPcondition2}
	For $b \in \R$, we say that $\mc{I}$ satisfies the Palais-Smale-Pohozaev condition at level $b$ (shortly the $(PSP)_b$ condition), if for any sequence $(\lambda_n, u_n) \subset \R \times H^s_r(\R^N)$ such that 
	\begin{equation}\label{prima}
	\mc{I} (\lambda_n, u_n) \to b,
	\end{equation}
	\begin{equation}\label{seconda}	
\partial_{\lambda} 	\mc{I}(\lambda_n, u_n) \to 0, 
	\end{equation}
	\begin{equation}\label{seconda3}
\norm{\partial_u 	\mc{I}(\lambda_n, u_n)}_{(H^s_r(\R^N))^*} \to 0,
\end{equation}
	\begin{equation}\label{terza4.14}
	\mc{P}(\lambda_n, u_n) \to 0,
	\end{equation}
	it happens that $(\lambda_n, u_n)$ has a strongly convergent subsequence in $\R \times H^s_r(\R^N)$.
\end{definition}

We will show the following result.

\begin{proposition}\label{PSP}
	Assume 
	\textnormal{(g1)--(g3)}.	Let $b \in \R$, $b <0$. Then $\mc{I}$ satisfies the $(PSP)_b$ condition on $\R\times H_r^s(\R^N)$.
\end{proposition}

\claim Proof.
Let $b \in \R$, $b <0$ and suppose that $(\lambda_n, u_n) \subset \R \times H^s_r(\R^N)$ satisfies \eqref{prima}--\eqref{terza4.14}. 
We will show that $(\lambda_n, u_n)$ has a strongly convergent subsequence in several steps.

\textbf{Step 1:} $\lambda_n$ is bounded from below. Indeed
\begin{eqnarray*}
	\frac{m}{2} e^{\lambda_n} &=&\frac{1}{N} \mc{P}(\lambda_n, u_n) - \mc{I}(\lambda_n, u_n) + \frac{s}{N}\norm{(-\Delta)^{s/2}u_n }_2^2 \\
	&\geq& \frac{1}{N} \mc{P}(\lambda_n, u_n) -\mc{I}(\lambda_n, u_n)
\end{eqnarray*}
hence
$$\frac{m}{2} \liminf_n e^{\lambda_n} \geq 0 -b >0,$$
which implies (since $m>0$) that $\lambda_n$ is bounded from below. 

\textbf{Step 2:} $\norm{u_n}_2^2\to m$. Indeed we have
$$\partial_{\lambda} \mc{I}(\lambda_n, u_n) = \frac{e^{\lambda_n}}{2} \left(\norm{u_n}_2^2 -m\right)\to 0,$$
which implies the claim by Step 1.

\textbf{Step 3:} $\norm{(-\Delta)^{s/2}u_n}_2^2$ and $\lambda_n$ are bounded (from above) as $n\to +\infty$. Indeed, by \eqref{maggiorazione-0} and the fractional Gagliardo-Nirenberg inequality \cite{Par0} we have 
\begin{eqnarray*}
	\abs{\partial_{u}\mc{I}(\lambda_n, u_n)u_n} &\geq& \norm{(-\Delta)^{s/2}u_n}_2^2 + e^{\lambda_n} \norm{u_n}_2^2 - \intRN |g(u_n)u_n| \\
	&\geq& \norm{(-\Delta)^{s/2}u_n}_2^2 + \left(e^{\lambda_n} - C_{\delta}\right) \norm{u_n}_2^2 - \delta \norm{u_n}_{p+1}^{p+1} \\
	&\geq& \norm{(-\Delta)^{s/2}u_n}_2^2 + \left(e^{\lambda_n} - C_{\delta}\right) \norm{u_n}_2^2 - \delta C \norm{(-\Delta)^{s/2}u_n}_2^2 \norm{u_n}_2^{p-1};
\end{eqnarray*}
moreover
\begin{eqnarray*}
	\abs{\partial_{u} \mc{I}(\lambda_n, u_n)u_n} &\leq& \norm{\partial_{u}I(\lambda_n, u_n)}_{(H^s_r(\R^N))^*} \norm{u_n}_{H^s_r(\R^N)}\\
	&=& \norm{\partial_{u}I(\lambda_n, u_n)}_{(H^s_r(\R^N))^*} \sqrt{\norm{(-\Delta)^{s/2} u_n}_2^2 + \norm{u_n}_2^2}.
\end{eqnarray*}
Set $\varepsilon_n= \norm{\partial_{u}I(\lambda_n, u_n)}_{(H^s_r(\R^N))^*}$ and (by Step 2) $\norm{u_n}_2^2= m + o(1)$, we finally have, joining the previous two inequalities, that
\begin{eqnarray*}
\lefteqn{\norm{(-\Delta)^{s/2} u_n}_2^2 \left(1-\delta C(m+o(1))^{\frac{p-1}{2}}\right) + \left(e^{\lambda_n} - C_{\delta}\right) ( m + o(1))}\\ 
&&\leq \varepsilon_n \sqrt{\norm{(-\Delta)^{s/2} u_n}_2^2 + m + o(1)}.
\end{eqnarray*}
Choosing $\delta>0$ small so that $\delta Cm^{\frac{p-1}{2}}<1$, we obtain the claim.

\textbf{Step 4:} Conclusion. By Steps 1-3, we have that $(\lambda_n, u_n)$ is bounded in $\R\times H^s_r(\R^N)$. 
Hence, up to a subsequence, $\lambda_n \to \lambda$ and $u_n \rightharpoonup u$ in $H^s_r(\R^N)$. 
Therefore, we obtain
$$\intRN g(u_n)u_n \to \intRN g(u)u \quad \hbox{ and } \quad \intRN g(u_n)u \to \intRN g(u)u.$$ 
Again by the assumption $\partial_{u}\mc{I}(\lambda_n, u_n)\to 0$ we obtain
\begin{eqnarray}
	0&=&\lim_n \partial_{u}\mc{I}(\lambda_n, u_n)u \notag \\ 
	&=& \lim_n \left(\intRN (-\Delta)^{s/2}u_n (-\Delta)^{s/2}u +e^{\lambda_n} \intRN u_n u - \intRN g(u_n)u\right) \notag \\
	&=& \norm{(-\Delta)^{s/2} u}_2^2 + e^{\lambda} \norm{u}_2^2 - \intRN g(u)u. \label{eq_dim_PSP1}
\end{eqnarray}
Since $\partial_{u}\mc{I} (\lambda_n, u_n)\to 0$ and $u_n \rightharpoonup u$, we have $\partial_{u}\mc{I} (\lambda_n, u_n)u_n\to 0$; thus
\begin{eqnarray}
	0&=&\lim_n \partial_{u}\mc{I}(\lambda_n, u_n)u_n \notag \\
	 &=& \lim_n \left(\norm{(-\Delta)^{s/2}u_n}_2^2 +e^{\lambda_n} \norm{u_n}_2^2- \intRN g(u_n)u_n\right) \notag \\
	&=& \lim_n\left(\norm{(-\Delta)^{s/2} u_n}_2^2 + e^{\lambda_n} \norm{u_n}_2^2\right) - \intRN g(u)u \label{eq_dim_PSP2}
\end{eqnarray}
and hence, joining \eqref{eq_dim_PSP1} and \eqref{eq_dim_PSP2},
$$\norm{(-\Delta)^{s/2} u_n}_2^2 + e^{\lambda_n} \norm{u_n}_2^2 \to \norm{(-\Delta)^{s/2} u}_2^2 + e^{\lambda} \norm{u}_2^2, $$
which easily implies (since $e^{\lambda_n}\to e^{\lambda}$ and $\norm{u_n}_2^2$ is bounded)
$$\norm{u_n}_{\lambda}^2 \to \norm{u}_{\lambda}^2, $$
where $\norm{\cdot}_{\lambda}^2= \norm{(-\Delta)^{s/2} \cdot}_2 + e^{\lambda} \norm{\cdot}_2^2$ is an equivalent norm on $H^s_r(\R^N)$. 
This, together with $u_n \rightharpoonup u$ in $H^s_r(\R^N)$ and the fact that $H^s_r(\R^N)$ is a Hilbert space, gives $u_n\to u$ strongly in $H^s_r(\R^N)$. 
\QED

\medskip


\begin{corollary}\label{PSP-cor}
	Assume \textnormal{(g1)--(g3)}. 
	Let $b \in \R$, $b <0$. Then $K^{PSP}_b \cap (\R \times \{0\}) = \emptyset$ and $K^{PSP}_b$ is compact.
\end{corollary}
\claim Proof.
Since 	$\partial_\lambda \mc{I}(\lambda, 0) = - \frac{e^\lambda}{2m} \neq 0$, we have $K^{PSP}_b \cap (\R \times \{0\}) = \emptyset$.
Proposition \ref{PSP} implies that $K^{PSP}_b$ is compact. 
\QED

\medskip

\begin{remark}
	We emphasize that the $(PSP)_b$ condition does not hold at level $b=0$. Indeed we can consider the unbounded sequence $(\lambda_j, 0)$ with $\lambda_j \to - \infty$ such that
	$$ \mc{I}(\lambda_j, 0) = \partial_\lambda \mc{I}(\lambda_j, 0) = - \frac{e^{\lambda_j}}{2} m \to 0$$ 
	and 
	$$ \partial_u \mc{I}(\lambda_j, 0) = 0, \quad \mc{P}(\lambda_j, 0) =0. $$ 
\end{remark}


\subsection{An augmented functional}

Following \cite{Jea0,HIT,HT} we introduce the augmented functional $\mc{H}: \R \times \R \times H^s_r(\R^N) \to \R$
\begin{equation}\label{ugual}
\mc{H}(\theta, \lambda, u)= \mc{I}(\lambda, u(e^{-\theta}\cdot)) .
\end{equation}
By the scaling properties of $\mc{I}$ we can recognize that
\begin{equation}\label{eq:22}
\mathcal{H}(\theta, \lambda, u)= 
\frac{e^{(N-2s)\theta}}{2}
\int_{\R^N} |(-\Delta)^{s/2}u|^2 - e^{N \theta}
\int_{\R^N} G(u)\, + 
\frac{e^\lambda}{2} \bigl( e^{N \theta}\|u\|_2^2 -m \bigr) 
\end{equation}
for all $(\theta, \lambda,u ) \in \R \times \R \times H^s_r(\R^N).$

Moreover, by standard calculations we have the following proposition.
\begin{proposition}
	For all $(\theta, \lambda,u ) \in \R \times \R \times H^s_r(\R^N)$, $h \in H^s_r(\R^N)$, $\beta \in \R$, we have
	\begin{itemize}
	\item[\textnormal{(i)}]
		$\partial_\theta \mc{H}(\theta, \lambda, u) = \mc{P}(\lambda, u(x /e^\theta)),	$		
	\item[\textnormal{(ii)}]			
		$\partial_\lambda \mc{H}(\theta, \lambda, u) = \partial_\lambda \mc{I}(\lambda, u(x/ e^\theta)), $
	\item[\textnormal{(iii)}]		 
		$\partial_u \mc{H}(\theta, \lambda, u) h(x) = \partial_u \mc{I}(\lambda, u(x/ e^\theta)) h(x/e^\theta), $		
	\item[\textnormal{(iv)}]
		$ \mc{H}(\theta + \beta, \lambda, u(e^\beta x)) = \mc{H}(\theta, \lambda, u). $		
	\end{itemize}
\end{proposition}

Now we define a metric on the Hilbert manifold
$$M= \R \times \R \times H^s_r(\R^N)$$
by setting 
\begin{eqnarray*}
{\|(\alpha, \nu, h)\|}_{(\theta,\lambda, u)}^2 &=&\abs{\left(\alpha, \nu,\norm{h(e^{-\theta} \cdot)}_{H^s_r(\R^N)}\right)}^2 \\
&=& \alpha^2 + \nu^2 + e^{N\theta} \norm{h}_2^2+ e^{(N-2s)\theta} \norm{(-\Delta)^{s/2} h}_2^2 
\end{eqnarray*}
for any $(\alpha, \nu, h) \in T_{(\theta,\lambda,u)} M = \R \times \R \times H^s_r(\R^N)$.
We also denote the dual norm on 
$T^*_{(\theta,\lambda,u)}M$ by $\|\cdot \|_{(\theta,\lambda, u), *}$. 
We notice that 
${\|(\cdot, \cdot, \cdot)\|}_{(\theta,\lambda, u)}^2$ depends only on $\theta$ and we can write ${\|(\cdot, \cdot, \cdot )\|}_{(\theta,\cdot, \cdot)}^2$.
Moreover for any $(\alpha, \nu, h) \in T_{(\theta,\lambda,u)}M$ and $\beta \in \R$ we have
\begin{equation}\label{eq_shift_norma}
{\|(\alpha, \nu, h(e^\beta x))\|}_{(\theta + \beta,\cdot, \cdot)}^2=
{\|(\alpha, \nu, h)\|}_{(\theta,\cdot, \cdot)}^2.
\end{equation}
Furthermore we define the standard distance between two points as the infimum of length of curves connecting the two points, namely
$$ \dist_M\big((\theta_0, \lambda_0, h_0), (\theta_1, \lambda_1, h_1)\big)= 
\inf_{\gamma \in \mathcal{G}} \int_0^1 \|\dot \gamma(t)\|_{\gamma(t)} dt $$
where 
$\mathcal{G}=\left\{\gamma \in C^1([0,1],M) \, \middle | \, \gamma(0)= (\theta_0, \lambda_0, h_0), \gamma(1)= (\theta_1, \lambda_1, h_1) \right\}.$

Observe that, if $\sigma$ is a path connecting $(\alpha_0, \nu_0, h_0)$ and $(\alpha_1, \nu_1, h_1)$, then by \eqref{eq_shift_norma} $\tilde{\sigma}(t)=(\sigma_1(t)+\beta, \sigma_2(t), (\sigma_3(t))(e^{\beta}\cdot))$ is a path connecting $(\alpha_0 +\beta, \nu_0, h_0(e^{\beta}\cdot))$ and $(\alpha_1+\beta, \nu_1, h_1(e^{\beta}\cdot))$ with same length, and hence
\begin{equation}\label{eq_shift_distance}
\dist_M\big((\alpha_0, \nu_0, h_0), (\alpha_1, \nu_1, h_1)\big) = \dist_M\big((\alpha_0 +\beta, \nu_0, h_0(e^{\beta}\cdot)), (\alpha_1+\beta, \nu_1, h_1(e^{\beta}\cdot))\big).
\end{equation}

Denote now $D=(\partial_\theta,\partial_\lambda,\partial_u)$ the gradient with respect to all the variables; a direct computation shows that
$$D\mathcal{H}(\theta, \lambda, u)(\alpha,\nu,h) 
= \mc{P} (\lambda, u(e^{-\theta} \cdot))\alpha +\partial_{\lambda} \mc{I}(\lambda, u(e^{-\theta}\cdot))\nu
+\partial_u \mathcal{I}(\lambda, u(e^{-\theta}\cdot))h(e^{-\theta} \cdot)$$
and thus we obtain
\begin{eqnarray*}
\lefteqn{\|{D\mathcal{H}(\theta, \lambda, u)\|}_{(\theta, \lambda, u),*}^2 }\\
	&=& \left|\left(\mathcal{P}(\lambda, u(e^{-\theta} \cdot)), \partial_{\lambda}\mathcal{I}(\lambda, u(e^{-\theta}\cdot)), \norm{\partial_u \mathcal{I}(\lambda, u(e^{-\theta}\cdot))}_{(H^s_r(\R^N))^*}\right)\right|^2 \\
	&=& \abs{\mathcal{P}(\lambda, u(e^{-\theta} \cdot))}^2 + \abs{\partial_{\lambda}\mathcal{I}(\lambda, u(e^{-\theta}\cdot))}^2 
		+ \norm{\partial_u \mathcal{I}(\lambda, u(e^{-\theta}\cdot))}_{(H^s_r(\R^N))^*}^2 .
\end{eqnarray*}
Now defined
$$\tilde{K}_b =\big \{ (\theta, \lambda, u) \in M \mid \mc{H}(\theta, \lambda, u)=b,\, D \mc{H}(\theta, \lambda,u)=0\big\}$$
the set of critical points at level $b$ of $\mathcal{H}$, we deduce
\begin{equation}\label{eq_confronto_K}
\tilde{K}_b = \big\{(\theta, \lambda, u(e^{\theta} \cdot)) \mid (\lambda, u)\in K^{PSP}_b, \; \theta \in \R\big\}.
\end{equation}

\begin{proposition}\label{PSPtilde}
	Assume \textnormal{(g1)--(g3)}. 
	Let $b \in \R$, $b <0$. Then the functional $\mathcal{H}$ satisfies the following Palais Smale type condition $(\widetilde{PSP})_b$. 
	That is, for each sequence $(\theta_n, \lambda_n, u_n)$ such that
	$$\mathcal{H}(\theta_n, \lambda_n, u_n) \to b,$$
	$$\norm{D \mathcal{H}(\theta_n, \lambda_n, u_n)}_{(\theta_n, \lambda_n, u_n),*} \to 0,$$
	we have, up to a subsequence,
	$$\dist_M((\theta_n, \lambda_n, u_n), \tilde{K}_b)\to 0.$$
\end{proposition}
We note that $(\widetilde{PSP})_b$ condition is different from the standard Palais-Smale condition and it ensures the compactness of $(\theta_n,\lambda_n,u_n)$ after a suitable scaling. 
We also highlight that, if $\tilde K_b\not=\emptyset$, then $\tilde K_b$ is not compact (see \eqref{eq_confronto_K}).

\smallskip

\claim Proof.	
Let $(\theta_n, \lambda_n, u_n)$ as in $(\widetilde{PSP})_b$. Then set $\hat{u}_n(x)= u_n(e^{-\theta_n}x)$ we have
$$ \mc{P}(\lambda_n, \hat{u}_n)\to 0,$$
$$\partial_{\lambda}\mathcal{I}(\lambda_n, \hat{u}_n) \to 0,$$
$$ \norm{\partial_u \mathcal{I}(\lambda_n,\hat{u}_n)}_{(H^s_r(\R^N))^*}\to 0,$$
and thus by Proposition \ref{PSP} the sequence $(\lambda_n, \hat{u}_n)$ is convergent (up to subsequences) to a $(\lambda, \hat{u})\in K^{PSP}_b$. 
Observe that, for each $n$, set $v_n(x)=\hat{u}(e^{\theta_n}x)$, we have $(\theta_n, \lambda, v_n)\in \tilde{K}_b$. 
Therefore by \eqref{eq_shift_distance}
\begin{eqnarray*}
	\dist_M((\theta_n, \lambda_n, u_n), \tilde{K}_b) &\leq& \dist_M((\theta_n, \lambda_n, u_n), (\theta_n, \lambda, v_n)) \\
	&=& \dist_M((0,\lambda_n, \hat{u}_n), (0, \lambda, \hat{u})) \\
	&\leq& \sqrt{|\lambda_n -\lambda|^2 + \norm{\hat{u}_n-\hat{u}}_{H^s_r(\R^N)}^2} \to 0,
\end{eqnarray*}
which reaches the claim.
\QED

\bigskip

\noindent
{\bf Notation.}
We use the following notation: for $\tilde A\subset M$ and $\rho>0$ we set
	$$ \tilde N_\rho(\tilde A) = \{ (\theta,\lambda,u)\in M\,|\, \dist_M((\theta,\lambda,u),\tilde A) <\rho\}, $$
while for $A\subset\R\times H_r^s(\R^N)$ and $R>0$ we set
	$$ N_R(A) = \{(\lambda,u)\in \R\times H_r^s(\R^N)\,|\, d((\lambda,u),A)<R\}, $$
where
	$$ d((\lambda,u),(\lambda',u')) = (|\lambda-\lambda'|^2 +\norm{u-u'}_{H_r^s}^2)^{1/2}. $$
We also write for $-\infty<a<b<\infty$
\begin{eqnarray*}
 \mc{I}^b &=& \{ (\lambda,u)\in\R\times H_r^s(\R^N)\,|\, \mc{I}(\lambda,u) \leq b\}, \\
 \mc{I}^b_a &=& \{ (\lambda,u)\in\R\times H_r^s(\R^N)\,|\, a\leq \mc{I}(\lambda,u) \leq b\}, \\
 \mc{H}^b &=& \{ (\theta,\lambda,u)\in M\,|\, \mc{H}(\theta,\lambda,u) \leq b\}, \\
 \mc{H}^b_a &=& \{ (\theta,\lambda,u)\in M\,|\, a\leq \mc{H}(\theta,\lambda,u) \leq b\}.
\end{eqnarray*}
Using these notation, as a corollary to Proposition \ref{PSPtilde}, we have
\begin{corollary}\label{cor_PStilde}
For any $\rho>0$ there exists a $\delta_\rho>0$ such that
		\begin{equation}\label{eq_curv_pend_gen}
		\forall \, (\theta,\lambda,u) \in \mc{H}^{b+\delta_{\rho}}_{b-\delta_{\rho}}\setminus \tilde N_\rho(\tilde K_b)
 \; : \; \norm{D\mc{H}(\theta,\lambda,u)}_{(\theta,\lambda,u),*}> \delta_{\rho}.
		\end{equation}
Here, if $\tilde{K}_b=\emptyset$, we regard $\tilde N_\rho(\tilde K_b)=\emptyset$.
\end{corollary}


\setcounter{equation}{0} 
\section{Construction of a deformation flow}\label{section:5}

Arguing as in Proposition 6.2 in \cite{HT} (see also \cite{IT}), we aim to prove the following Deformation Theorem in the fractional framework.

\begin{theorem}\label{defarg}
Let $b<0$, and assume $K_b^{PSP}= \emptyset$. 
Let $\bar{\varepsilon}>0$, then there exist $\varepsilon \in (0,\bar{\varepsilon})$ and $\eta: [0,1]\times (\R\times H^s_r(\R^N))\to \R\times H^s_r(\R^N)$ continuous such that
\begin{enumerate}
\item $\eta(0, \cdot,\cdot)=id_{\R\times H^s_r(\R^N)}$;
\item $\eta$ fixes $\mc{I}^{b-\bar{\varepsilon}}$, that is, $\eta(t, \cdot,\cdot)=id_{\mc{I}^{b-\bar{\varepsilon}}}$
for all $t\in [0,1]$;
\item $\mc{I}$ is non-increasing along $\eta$, and in particular $\mc{I}(\eta(t,\cdot, \cdot))\leq \mc{I}(\cdot, \cdot)$ for all $t \in [0,1]$;
\item $\eta(1, \mc{I}^{b+\varepsilon})\subseteq \mc{I}^{b-\varepsilon}$.
\end{enumerate}

\end{theorem}
We omit the proof of the Theorem since it will be very similar to the one made in the case of multiplicity (see Theorem \ref{thm_def_gen}). 
We remark that this deformation flow is not $C^1$ and it does not satisfy the two properties of the standard deformation flow, in general:
\begin{itemize}
	\item[(1)] $\eta(s+t, \lambda, u)= \eta(t, \eta(s,\lambda, u))$ with $s+t \in [0,1], (\lambda,u) \in \R \times H^s_r(\R^N)$;	
	\item[(2)] for $t \in [0,1]$, the map $(\lambda,u) \mapsto \eta(t,\lambda, u)$ is a homeomorphism.
\end{itemize}
We refer to Remark 3.2 in \cite{HT}.

We also stress that the deformation argument in Theorem \ref{defarg} works for $K^{PSP}_b$ but not for $K_b$ and thus, if $K^{PSP}_b=\emptyset$, then we have the statement (4) in Theorem \ref{defarg} even if $K_b\not=\emptyset$. 
We derive the following corollary (see also Corollary \ref{coroll_esist_Pm}).

\begin{corollary}\label{dedu}.
	Let $\bar b <0$ be a MP minimax value for $\mathcal{I}$. 
	Then $K^{PSP}_{\bar b} \not=\emptyset$, that is, $\mc{I}$ has a critical point $(\bar \lambda, \bar u)$ satisfying the Pohozaev identity, namely $\mc{P}(\bar \lambda, \bar u)=0$.
\end{corollary}


\setcounter{equation}{0} 
\section{Minimax critical points in the product space}\label{section:6}

For any $m >0$, let $B_m$ and $B_m'$ be the constants defined in Proposition \ref{negativo}, namely
$$ B_m = \inf_{\lambda<\lambda_0} \left(a(\lambda) -{e^\lambda\over 2}m \right), \quad B'_m=\inf_{(\lambda,u)\in\partial\Omega}\mc{I}(\lambda,u).$$
As a minimax class for $\mc{I}(\lambda,u)$, we define
\begin{eqnarray*}
	\Gamma^m = \big\{\xi\in C\big([0,1], \R \times H^s_r(\R^N)\big) & \mid 
	& \xi(0) \in \R \times \{ 0\},\ \mc{I}(\xi(0)) \leq B_m -1,\\ &&
	\xi(1) \not \in \Omega,		
	\ \mc{I}(\xi(1)) \leq B_m -1\big\}. 
\end{eqnarray*}

We have the following result.

\begin{proposition}\label{tom} Assume \textnormal{(g1)--(g3)}.
	\begin{itemize}
	\item[\textnormal{(i)}]
				 For any $m>0$, we have $\Gamma^m \neq \emptyset$.
	\item[\textnormal{(ii)}]	
		For sufficiently large $m>0$ there exists $\xi\in \Gamma^m$ such that
		\begin{equation}\label{terza}
		\max_{t\in [0,1]} \mc{I}(\xi(t)) <0.
		\end{equation}
	\item[\textnormal{(iii)}]	
		Assume \textnormal{(g4)}. 
		Then for any $m>0$ there exists $\xi\in \Gamma^m$ with the property \eqref{terza}.
	\end{itemize}
\end{proposition}

\claim Proof.
Let $\lambda_0\in (-\infty,\infty]$ be defined in \eqref{lambdazero}. 
For any $\lambda<\lambda_0$ we show there exists a path $\psi_\lambda\in\Gamma^m$ such that
\begin{equation} \label{aa}
	\max_{t\in [0,1]} \mc{I}(\psi_\lambda(t)) \leq a(\lambda)-{e^\lambda\over 2}m.
\end{equation}
Let $u_\lambda$ be a MP solution of $\partial_u\mc{J}(\lambda,u)=0$ (by Theorem \ref{Byeonminimizers}). 
Set $\zeta_\lambda(t)=u_\lambda(\cdot/t)$ for $t>0$ and $\zeta_\lambda(0)=0$ and note that, since $u_{\lambda}$ satisfies the Pohozaev identity, we have $\mc{I} (\lambda,\zeta_\lambda(t))
\to-\infty$ and $\mc{P}(\lambda,\zeta_\lambda(t)) \to-\infty$ as $t\to+\infty$. We can find $\gamma_\lambda=\zeta_\lambda(L\cdot)$ for $L\gg 1$ satisfying
$$ a(\lambda) = \max_{t\in [0,1]} \mc{J}(\lambda,\gamma_\lambda(t)),$$
$$ \mc{I}(\lambda,\gamma_\lambda(1)) \leq B_m-1, \quad \gamma_{\lambda}(1)\notin \Omega.$$
We also note that $t\mapsto \mc{I}(t,0)=-{e^t\over 2}m$ is decreasing and tending to $-\infty$
as $t\to+\infty$. Thus, joining $\gamma_\lambda$ and $t\mapsto (\lambda+Lt,0);\, [0,1]\to\R\times H_r^s(\R^N)$
for $L\gg 1$, we find a path $\psi_\lambda\in\Gamma^m$, defined as
$$	\psi_{\lambda}(t) = \parag{
		(\lambda+L(1-2t), \,0) & \quad \hbox{ if $t \in [0,1/2]$,} \\ 
		\left(\lambda, \, \gamma_{\lambda}(2t-1) \right) & \quad \hbox{ if $t \in (1/2,1]$}
		}
$$
 with \eqref{aa}.
Thus in particular we have (i).

Next we deal with (ii) and (iii). By \eqref{aa}, we have that (ii) follows easily; (iii) also follows from Proposition \ref{out}.
\QED

\bigskip

We notice that each path in $\Gamma^m$ passes through $\partial \Omega$, thus the minimax value
\begin{equation}\label{crti}
b_m = \inf_{\xi\in\Gamma^m}\max_{t\in [0,1]} \mc{I}(\xi(t))
\end{equation} 
verifies $b_m\geq B'_m$ and hence by Proposition \ref{negativo} it is well-defined and finite. 
Since Palais-Smale-Pohozaev condition holds on $(-\infty,0)$, it is important to estimate $b_m$. 
We have the following result.
\begin{proposition}\label{blue}
	 Assume \textnormal{(g1)--(g3)}.
	\begin{itemize}
		\item[\textnormal{(i)}]
		There exists $m_0>0$ such that
			$$ b_m<0 \quad \hbox{for $m>m_0$}. $$
		\item[\textnormal{(ii)}]Assume \textnormal{(g4)} in addition, then $m_0=0$, that is, 
			$$ b_m <0 \quad \hbox{for all $m>0$}.$$
	\item[\textnormal{(iii)}]	
		We have $b_m=B'_m=B_m$.
	\end{itemize}
\end{proposition}

\claim Proof.
By \eqref{aa}, we have
\begin{eqnarray}\label{mor}
 b_m &\leq& a(\lambda)-{e^\lambda\over 2}m \nonumber \\
 &=& e^\lambda\left({a(\lambda)\over e^\lambda}-{m\over 2}\right)
 \quad \hbox{for all}\ \lambda<\lambda_0.
\end{eqnarray}
Setting
$$ m_0 = 2 \inf_{\lambda<\lambda_0}{a(\lambda)\over e^\lambda}\geq 0, $$
we have $b_m<0$ for $m>m_0$. 
Thus we have (i). 
By Proposition \ref{out}, we have $m_0=0$ under the assumption \textnormal{(g4)} and thus we have (ii).

Finally, from \eqref{mor} it follows $b_m \leq B_m$. As already observed $b_m \geq B'_m \geq B_m$, from which we deduce (iii).
\QED

\bigskip

By Proposition \ref{blue} and Corollary \ref{dedu} we conclude that the level $b_m$, defined in \eqref{crti}, is a critical value of $\mc{I}$ in the product space $\R \times H^s_r(\R^N)$ and thus Theorem \ref{S:1.1} and Theorem \ref{S:1.12} hold.

\begin{corollary}\label{coroll_esist_Pm}
Let $m>m_0$.
 Then there exists a solution of problem \eqref{problem} which satisfies the Pohozaev identity.
If moreover \textnormal{(g4)} holds, then there exists a solution of \eqref{problem} for each $m>0$.
\end{corollary}

\claim Proof.
Let $\bar{\varepsilon}\in (0, 1)$. By Theorem \ref{defarg}, in correspondence to $b_m<0$, there exists $\varepsilon\in (0,\bar{\varepsilon})$ and $\eta$ satisfying $1)-4)$. By definition of inf, there exists $\gamma \in \Gamma^m$ such that
$$\max_{t\in [0,1]} \mc{I}(\gamma(t)) < b_m + \varepsilon,$$
that is 
\begin{equation}\label{d1}
\gamma([0,1])\subseteq \mc{I}^{b_m+\varepsilon}.
\end{equation}
Set
$$\tilde{\gamma}(t)=\eta(1,\gamma(t)),$$
we show that $\tilde{\gamma}\in \Gamma^m$. Indeed for $i\in\{0,1\}$, since $\mc{I}(\gamma(i))\leq B_m - 1 \leq b_m-\bar{\varepsilon}$, 
Theorem \ref{defarg} implies that $\tilde{\gamma}(i)=\eta(1,\gamma(i))=\gamma(i)\in \mc{I}^{b_m-\bar{\varepsilon}}$, and thus $\tilde{\gamma}(0)=\gamma(0)\in\R\times\{ 0\}$, $\tilde{\gamma}(1)=\gamma(1)\not\in\Omega$.
Therefore 
\begin{equation}\label{d2}
b_m \leq \max_{t\in [0,1]} \mc{I}(\tilde{\gamma}(t)).
\end{equation}
By contradiction, assume $K_{b_m}^{PSP} = \emptyset$. By the properties of $\eta$ and \eqref{d1} we obtain that $\tilde\gamma([0,1])=\eta(1,\gamma([0,1]))\subseteq \mc{I}^{b_m-\varepsilon}$, that is
$$\max_{t\in [0,1]} \mc{I}(\eta(1,\gamma(t)))\leq b_m-\varepsilon.$$
This is in contradiction with \eqref{d2}, and we conclude the proof.
\QED

\bigskip

In Theorems \ref{S:1.1} and \ref{S:1.12} we find solutions via mountain pass minimax methods. 
We remark that these solutions are characterized as minimizers of the functional $\mc{L}$ on $\mc{S}_m$, where $\mc{L}:\, H_r^s(\R^N)\to\R$ is defined by
	$$	\mc{L}(u)=\half\norm{(-\Delta)^{s/2}u}_2^2
		-\int_{\R^N} G(u)
	$$
and $\mc{S}_m$ is the $L^2$-sphere in $H_r^S(\R^N)$, i.e. 
$$\mc{S}_m=\{u\in H_r^s(\R^N)\,|\, \norm u_2^2=m\}.$$
Setting
	$$	\kappa_m=\inf_{u\in\mc{S}_m} \mc{L}(u),
	$$
we have the following result.
\begin{proposition} \label{minimizing}
Under the assumption of Theorem \ref{S:1.1}, we have for $m>m_0$,
\begin{itemize}
\item[\textnormal{(i)}] $-\infty<\kappa_m<0$ and $\kappa_m$ is attained.
\item[\textnormal{(ii)}] $\kappa_m=b_m$, where $b_m$ is defined in \eqref{crti}.
\end{itemize}
Moreover, in the assumption of Theorem \ref{S:1.12}, $m_0=0$.
\end{proposition}

\claim Proof.
\\
\textbf{Step 1:} $\kappa_m>-\infty$ and $\kappa_m<0$ for $m>m_0$.

\smallskip

\noindent
By (g1)-(g2), for any $\delta>0$ there exists $C_\delta>0$ such that
	$$	\mc{L}(u) \geq \half\norm{(-\Delta)^{s/2}u}_2^2
		-\frac\delta{p+1}\norm u_{p+1}^{p+1}
		-C_\delta \norm u_2^2.
	$$
By the fractional Gagliardo-Nirenberg inequality we have, for $u\in \mc{S}_m$
	\begin{eqnarray*}
	\mc{L}(u) &\geq \half\norm{(-\Delta)^{s/2}u}_2^2
		-\frac{C\delta}{p+1}\norm{(-\Delta)^{s/2}u}_2^2
			\norm u_2^{p-1}
		-C_\delta \norm u_2^2 \\
		&=\left(\half-\frac{C\delta}{p+1}m^{\frac{p-1}2}\right)
			\norm{(-\Delta)^{s/2}u}_2^2 -C_\delta m.
	\end{eqnarray*}
Choosing $\delta>0$ small so that $\half-\frac{C\delta}{p+1}m^{\frac{p-1}2}>0$, we have $\kappa_m\geq -C_\delta m>-\infty$.
\\ Since the solution $u_*\in\mc{S}_m$ obtained in Theorem \ref{S:1.1} satisfies
	$$	0>b_m=\mc{L}(u_*) \geq \kappa_m,
	$$
we have $\kappa_m<0$ for $m>m_0$.

\smallskip

\noindent
 \textbf{Step 2:} For $m>m_0$, $\kappa_m$ is attained.

\smallskip

\noindent
To show the existence of a minimizer of $\mc{L}$ on $\mc{S}_m$, we use a linear action $ \Phi:\, \R\to L(H_r^s(\R^N))$ defined by
	$$	(\Phi_\theta v)(x)= e^{\frac N2 \theta}v(e^\theta x).
	$$
We note that $\mc{S}_m$ is invariant under 
$\Phi_\theta$, that is, $\Phi_\theta(\mc{S}_m) = \mc{S}_m$. 
Let 
$$N=\R\times \mc{S}_m$$
and on the tangent bundle $TN=\R\times T\mc{S}_m =\coprod_{(\theta,u)\in N} (\R\times T_u\mc{S}_m)$
we introduce a $C^2$-metric
	$$	\norm{(\kappa,v)}_{(\theta,u)} 
	=\left( \kappa^2+\norm{\Phi_\theta v}_{H^s(\R^N)}^2
		\right)^{1/2}
	$$
for all $(\theta, u) \in N$ and $(\kappa, v) \in TN$. We also introduce $\widehat{\mc{L}}:\, N \to \R$ by
	$$	\widehat{\mc{L}}(\theta,u)
		= \mc{L}(\Phi_\theta u)
		= \half e^{2s\theta}\norm{(-\Delta)^{s/2}u}_2^2
		-e^{-N\theta}\int_{\R^N} G(e^{\frac N2\theta} u).
	$$
We note that 
	$$	\inf_{(\theta,u)\in N}\widehat{\mc{L}}(\theta,u)=\kappa_m.
	$$
Since $\kappa_m \in \R$ by Step 1, applying Ekeland's principle, there exists a sequence $(\theta_j,u_j)_{j=1}^\infty\subset N$ such that 
	\begin{eqnarray*}
	&\widehat{\mc{L}}(\theta_j,u_j)\to \kappa_m, \\
	&\norm{D\widehat{\mc{L}}(\theta_j,u_j)}_{T^*_{(\theta_j,u_j)}N}\to 0. 
	\end{eqnarray*}
That is, noting $T_u\mc{S}_m=\big\{ v\in H_r^s(\R^N) \mid \int_{\R^N} uv =0\big \}$,
	\begin{eqnarray*}
	&\partial_\theta\widehat{\mc{L}}(\theta_j,u_j)\to 0,\\
	&\norm{\partial_u \widehat{\mc{L}}
		(\theta_j,u_j)}_{T^*_{u_j}\mc{S}_m}
	=\sup_{
	v\in T_{u_j}\mc{S}_m, \,
			\norm{\Phi_{\theta_j}v}_{H^s(\R^N)} \leq 1 }
	\abs{\partial_u\widehat{\mc{L}}(\theta_j,u_j)v}
	\to 0.	
	\end{eqnarray*}
Setting $\widehat u_j=\Phi_{\theta_j}u_j$, we thus have
	\begin{eqnarray}
	&\norm{\widehat u_j}_2^2=m, \label{b.0}\\
	&\mc{L}(\widehat u_j)
	=\half\norm{(-\Delta)^{s/2}\widehat u_j}_2^2 
		-\int_{\R^N} G(\widehat u_j) \to \kappa_m, \label{b.1}\\
	&s\norm{(-\Delta)^{s/2}\widehat u_j}_2^2
		+N\int_{\R^N} G(\widehat u_j)
		-\frac N2 \int_{\R^N} g(\widehat u_j)\widehat u_j
	\to 0	\label{b.2}
 \end{eqnarray}
and for a suitable $\mu_j\in\R$
 \begin{equation}\label{b.3}
 \mc{L}'(\widehat u_j) \widehat{v}
 +\mu_j\int_{\R^N} \widehat u_j \widehat{v}= o(1)\norm{\widehat{v}}_{H^s(\R^N)}
 \quad\hbox{for all}\ \widehat{v}\in H_r^s(\R^N).
 \end{equation}
By using \eqref{b.1} and arguing as in Step 1 we see that $\widehat{u}_j$ is bounded in $H^s_r(\R^N)$. Thus, choosing $\widehat{v}= \widehat{u}_j$ in \eqref{b.3}, we have
	$$	\norm{(-\Delta)^{s/2}\widehat u_j}_2^2
	-\int_{\R^N} g(\widehat u_j)\widehat u_j + \mu_j m
	=o(1),
	$$
which, joined to \eqref{b.2}, gives
	$$	\left(1-\frac{2s}N\right)\norm{(-\Delta)^{s/2}\widehat u_j}_2^2
		-2\int_{\R^N} G(\widehat u_j) +\mu_j m =o(1),
	$$
that is
	$$	-\frac{2s}N \norm{(-\Delta)^{s/2}\widehat u_j}_2^2
		+ 2\mc{L}(\widehat u_j)
		+\mu_j m =o(1).
	$$
Thus, by \eqref{b.1},
	$$	\mu_j \geq -\frac{2\kappa_m}m + o(1)
	$$
which implies, by Step 1, that $\mu_j>0$ for $j$ large and hence we can write $\mu_j=e^{\widehat\lambda_j}$ for some $\widehat\lambda_j\in\R$. 
Relations \eqref{b.0}--\eqref{b.3} imply that $(\widehat\lambda_j,\widehat u_j)$ satisfies \eqref{prima}--\eqref{terza4.14} with $b=\kappa_m<0$.
Thanks to Proposition \ref{PSP}, $(\widehat \lambda_j,\widehat u_j)$ has a strongly convergent subsequence to some $(\widehat\lambda_*,\widehat u_*) \in N$, which shows the existence of a minimizer $\widehat{u}_*$. 
Thus (i) is proved.

\smallskip

\noindent
\textbf{Step 3:} For $m>m_0$, $\kappa_m=b_m$.

\smallskip

\noindent
In Step 1, we showed $b_m\geq \kappa_m$. On the other hand by the argument in Step 2, for the minimizer $\widehat u_*$ of $\mc{L}$ on $\mc{S}_m$, there exists $\widehat \lambda_*\in\R$ such that
	\begin{eqnarray*}
	&\mc{I}(\widehat\lambda_*,\widehat u_*)=\kappa_m,
	\quad 
	\partial_u\mc{I}(\widehat\lambda_*,\widehat u_*)=0,
		\\
	&\partial_\lambda\mc{I}(\widehat\lambda_*,\widehat u_*)=0,
	\quad 
	\mc{P}(\widehat\lambda_*,\widehat u_*)=0.
	\end{eqnarray*}
Thus, by the argument in previous sections, there exists $\xi_*\in\Gamma^m$ such that
	$$	\max_{t\in [0,1]}\mc{I}(\xi_*(t))
		= \mc{I}(\widehat \lambda_*,\widehat u_*)
		= \kappa_m,
	$$
which implies $b_m=\kappa_m$ and the proof of Proposition \ref{minimizing} is completed. 
\QED


\setcounter{equation}{0} 
\section{Multiple normalized solutions}\label{section:7}

In the whole section we assume, in addition, (g5).


\subsection{Deformation theorems}

In what follows we will use the following terminology. Set $\G=\Z_2$, we consider the action $\sigma$ of $\G$ on $\R\times H^s_r(\R^N)$ and on $M$, that is
$$\sigma:\,(\pm1,\lambda, u)\mapsto (\lambda, \pm u);\, \G\times (\R\times H^s_r(\R^N)) \to \R\times H^s_r(\R^N),$$
$$\sigma:\,(\pm1,\theta, \lambda, u)\mapsto (\theta, \lambda, \pm u);\, \G\times M \to M.$$
We notice that $\mc{I}$ and $\mc{H}$ are invariant under this action (i.e. they are even in $u$), as well as the set $\Omega$ (i.e. it is symmetric with respect the axis $\R$). 
We highlight instead that the function $\eta=(\eta_1, \eta_2):\R\times H^s_r(\R^N) \to \R\times H^s_r(\R^N)$ (resp. $\tilde{\eta}=(\tilde{\eta}_0,\tilde{\eta}_1,\tilde{\eta}_2):M \to M$) is equivariant if $\eta_1$ is even and $\eta_2$ is odd (resp. $\tilde{\eta}_0$ and $\tilde{\eta}_1$ are even and $\tilde{\eta}_2$ is odd).
We recall that a function $f$ is said to be invariant under the action $g \cdot x$ if $f(g \cdot x) = f(x)$, while it is equivariant if $f(g \cdot x) = g \cdot f(x)$.
We want to prove the following.


\begin{theorem}\label{thm_def_gen}
Let $b<0$, and let $\mc{O}$ be a neighbourhood of $K_b^{PSP}$. 
Let $\bar{\varepsilon}>0$, then there exist $\varepsilon \in (0,\bar{\varepsilon})$ and $\eta: [0,1]\times (\R\times H^s_r(\R^N))\to (\R\times H^s_r(\R^N))$ continuous such that
\begin{enumerate}
\item $\eta(0, \cdot,\cdot)=id_{\R\times H^s_r(\R^N)}$;
\item $\eta$ fixes $\mc{I}^{b-\bar{\varepsilon}}$, that is, $\eta(t, \cdot,\cdot)=id_{\mc{I}^{b-\bar{\varepsilon}}}$ for all $t \in [0,1]$;
\item $\mc{I}$ is non-increasing along $\eta$, and in particular $\mc{I}(\eta(t,\cdot, \cdot))\leq \mc{I}(\cdot, \cdot)$ for all $t \in [0,1]$;
\item if $K_b^{PSP}= \emptyset$, then $\eta(1, \mc{I}^{b+\varepsilon})\subseteq \mc{I}^{b-\varepsilon}$;
\item if $K_b^{PSP}\not=\emptyset$, then
$$\eta(1,\mc{I}^{b+\varepsilon}\setminus \mc{O}) \subseteq \mc{I}^{b-\varepsilon}$$
and
$$\eta(1,\mc{I}^{b+\varepsilon}) \subseteq \mc{I}^{b-\varepsilon}\cup \mc{O};$$
\item $\eta(t, \cdot,\cdot)$ is $\G$-equivariant, in the sense mentioned before.
\end{enumerate}
\end{theorem}

To prove this, we work first on the functional $\mc{H}$, for which we obtained a $(\widetilde{PSP})$ condition.

\begin{theorem}\label{thm_def_H}
Let $b<0$. 
Let moreover $\rho>0$ and write $\tilde{\mc{O}}=\tilde{N}_{\rho}(\tilde{K}_b)$. 
Let $\bar{\varepsilon}>0$, then there exist $\varepsilon \in (0,\bar{\varepsilon})$ and $\tilde{\eta}: [0,1]\times M\to M$ continuous such that
\begin{enumerate}
\item $\tilde{\eta}(0, \cdot,\cdot)=id_M$;
\item $\tilde{\eta}$ fixes $\mc{H}^{b-\bar{\varepsilon}}$, that is $\tilde{\eta}(t, \cdot,\cdot)=id_{\mc{H}^{b-\bar{\varepsilon}}}$ for all $t\in [0,1]$;
\item $\mc{H}$ is non-increasing along $\tilde{\eta}$, and in particular $\mc{H}(\tilde{\eta}(t,\cdot,\cdot, \cdot))\leq \mc{H}(\cdot, \cdot, \cdot)$ for all $t\in [0,1]$; 
\item if $\tilde{K}_b= \emptyset$, then $\tilde{\eta}(1, \mc{H}^{b+\varepsilon})\subseteq \mc{H}^{b-\varepsilon}$;
\item if $\tilde{K}_b \neq \emptyset$, then
$$\tilde{\eta}(1,\mc{H}^{b+\varepsilon}\setminus \tilde{\mc{O}}) \subseteq \mc{H}^{b-\varepsilon}$$
and
$$\tilde{\eta}(1,\mc{H}^{b+\varepsilon}) \subseteq \mc{H}^{b-\varepsilon}\cup \tilde{\mc{O}};$$
\item $\tilde{\eta}(t, \cdot,\cdot)$ is $\G$-equivariant, in the sense mentioned before.
\end{enumerate}
\end{theorem}

We postpone the proof of Theorem \ref{thm_def_H} for $\mc{H}$ and see now how to use it to deduce the one for $\mc{I}$. Introduce first the following notation:
$$\pi: M \to \R\times H^s_r(\R^N), \; \pi(\theta, \lambda, u) = (\lambda, u(e^{-\theta}\cdot)),$$
$$\iota: \R\times H^s_r(\R^N) \to M, \; \iota(\lambda, u) = (0,\lambda, u),$$
which are a kind of rescaling projection and immersion. Observe that
$$\pi \circ \iota = id_{\R\times H^s_r(\R^N)}, \quad \hbox{(while $\iota \circ \pi \neq id_M$),}$$
$$\mc{H}\circ \iota=\mc{I}, \quad \mc{I}\circ \pi=\mc{H},$$
$$\pi(\tilde{K}_b) = K_b^{PSP}.$$
For $\tilde{\eta}$ obtained in Theorem \ref{thm_def_H}, define "$\eta=\pi \circ \tilde{\eta} \circ \iota$" up to the time; more precisely
\begin{equation}\label{eq_def_flow}
\eta(t,\lambda, u)=\pi(\tilde{\eta}(t,\iota(\lambda, u))).
\end{equation}
It is now a straightforward computation showing that $\eta$ satisfies the requests of Theorem \ref{thm_def_gen}. 
A delicate issue, anyway, is to show the intuitive fact that neighbourhoods of $\tilde{K}_b$ are brought to neighbourhoods of $K_b^{PSP}$. More precisely we have the following result.

\begin{lemma}\label{lem_intorni}
Assume that $K_b^{PSP}$ is compact (for instance, $b<0$). Let $\rho>0$, then there exists $R(\rho)>0$ such that, set $\tilde{\mc{O}}=\tilde{N}_{\rho}(\tilde{K}_b)$ and $\mc{O}=N_{R(\rho)}(K_b^{PSP})$, we have
$$\pi(\tilde{\mc{O}}) \subset \mc{O},$$
i.e.
$$\dist_M((\theta, \lambda, u), \tilde{K}_b) \leq \rho \implies d((\lambda, u(e^{-\theta}\cdot)), K_b^{PSP}) \leq R(\rho).$$
In particular, for $\theta =0$ we have
\begin{equation}\label{eq_theta0}
\dist_M((0, \lambda, u), \tilde{K}_b) \leq \rho \implies d((\lambda, u), K_b^{PSP}) \leq R(\rho),
\end{equation}
that is
$$\iota(\complement \mc{O}) \subseteq \complement \tilde{\mc{O}}$$
where $\complement$ denotes the complement of a set. Moreover
$$\lim_{\rho\to 0} R(\rho)=0.$$
\end{lemma} 

\claim Proof.
We observe that is sufficient to prove \eqref{eq_theta0} since by \eqref{eq_shift_distance}
$$\dist_M((\theta,\lambda, u), \tilde{K}_b) = \dist_M((0,\lambda, u(e^{-\theta}\cdot)), \tilde{K}_b).$$
Let $\varepsilon>0$. By definition of $\dist_M((0,\lambda,u), \tilde K_b)$ there exists a $\sigma=\sigma(t)$, $\sigma=(\theta, \lambda, u)$, such that $\sigma(0)=(0,\lambda, u)$, $\sigma(1)\in \tilde{K}_b$ and
\begin{equation}\label{eq_dim_def_1}
\int_0^1 \norm{\dot{\sigma}(t)}_{\sigma(t)} dt \leq \rho + \varepsilon.
\end{equation}
By \eqref{eq_confronto_K} we have $(\lambda(1), u(1)(e^{-\theta(1)} \cdot))\in K_b^{PSP}$ and thus
\begin{eqnarray*}
\lefteqn{ \dist((\lambda, u), K_b^{PSP})}
\\ &\leq& \norm{(\lambda, u) - (\lambda(1), u(1)(e^{-\theta(1)} \cdot))}_{\R\times H^s_r(\R^N)}\\
&\leq & \norm{(\lambda, u) - (\lambda(1), u(1))}_{\R\times H^s_r(\R^N)}+ \norm{(\lambda(1), u(1)) - (\lambda(1), u(1)(e^{-\theta(1)} \cdot))}_{\R\times H^s_r(\R^N)}\\
&= & \norm{(\lambda(0), u(0)) - (\lambda(1), u(1))}_{\R\times H^s_r(\R^N)} + \norm{u(1) - u(1)(e^{-\theta(1)} \cdot)}_{H^s_r(\R^N)}\\
&=& I + II.
\end{eqnarray*}
Focus on $I$. We have, by the fundamental theorem of calculus and H\"older inequality,
\begin{eqnarray*}
I&=&\norm{(\lambda(0), u(0)) - (\lambda(1), u(1))}_{\R\times H^s_r(\R^N)} \leq \int_0^1 \left(\abs{\dot{\lambda}(t)}^2 + \norm{\dot{u}(t)}_{H^s_r(\R^N)}^2\right)^{1/2} dt \\
&= & \int_0^1 \left(\abs{\dot{\lambda}(t)}^2 + \norm{\dot{u}(t)}_{2}^2 + \norm{(-\Delta)^{s/2}\dot{u}(t)}_{2}^2\right)^{1/2} dt. 
\end{eqnarray*}
In order to use \eqref{eq_dim_def_1} it must appear the norm associated to $M$, which we recall is
$$ \norm{\dot{\sigma}(t)}_{\sigma(t)}^2= \dot{\theta}(t)^2 + \dot{\lambda}(t)^2 + e^{N\theta(t)} \norm{\dot{u}(t)}_2^2 + e^{(N-2s)\theta(t)} \norm{(-\Delta)^{s/2} \dot{u}(t)}_2^2.$$
Since we do not know the sign of $N\theta(t)$, 
we need an estimate on $\theta(t)$ and a corrective factor. Indeed, recalled that $\theta(0)=0$, we have
$$|\theta(t)| = |\theta(t)-\theta(0)| \leq \int_0^1 |\dot{\theta}(t)| dt \leq \int_0^1 \norm{\dot{\sigma}(t)}_{\sigma(t)} dt \leq \rho + \varepsilon.$$
Thus $\theta(t) \geq -(\rho + \varepsilon) \geq -\frac{N}{N-2s} (\rho + \varepsilon)$ which imply
$$e^{N(\rho+\varepsilon)} \geq 1, \quad e^{N(\rho+\varepsilon)} e^{N\theta(t)}\geq 1, \quad e^{N(\rho+\varepsilon)} e^{(N-2s)\theta(t)}\geq 1 $$
and hence we obtain
\begin{eqnarray*}
I& \leq& e^{\frac{N(\rho+\varepsilon)}{2}} \int_0^1 \left(\abs{\dot{\lambda}(t)}^2 + e^{N\theta(t)} \norm{\dot{u}(t)}_{2}^2 + e^{(N-2s)\theta(t)} \norm{(-\Delta)^{s/2}\dot{u}(t)}_{2}^2\right)^{1/2} dt \\
& \leq& e^{\frac{N(\rho+\varepsilon)}{2}} \int_0^1 \left(\abs{\dot{\theta}(t)^2} + \abs{\dot{\lambda}(t)}^2 + e^{N\theta(t)} \norm{\dot{u}(t)}_{2}^2 + e^{(N-2s)\theta(t)} \norm{(-\Delta)^{s/2}\dot{u}(t)}_{2}^2\right)^{1/2} dt \\
&=& e^{\frac{N(\rho+\varepsilon)}{2}} \int_0^1 \norm{\dot{\sigma}(t)}_{\sigma(t)} dt 
\leq e^{\frac{N(\rho+\varepsilon)}{2}} (\rho + \varepsilon) \stackrel{\varepsilon \to 0} \to e^{\frac{N\rho}{2}} \rho.
\end{eqnarray*}
Focus now on $II$. Set $\bar{\omega}= u(1)(e^{-\theta(1)}\cdot)$ we have $\bar{\omega} \in P_2(K_b^{PSP})$ (where $P_2$ is the projection on the second component) with $|\theta(1)|\leq \rho + \varepsilon$, and thus
\begin{eqnarray*}
II &=& \norm{u(1) - u(1)(e^{-\theta(1)} \cdot)}_{H^s_r(\R^N)} = \norm{\bar{\omega}(e^{\theta(1)}\cdot) - \bar{\omega}}_{H^s_r(\R^N)}\\
&\leq& \sup \left \{ \norm{\omega(e^{\alpha}\cdot)- \omega}_{H^s_r(\R^N)} \mid |\alpha| \leq \rho + \varepsilon, \; \omega \in P_2(K_b^{PSP}) \right \}.
\end{eqnarray*}
Since $P_2(K_b^{PSP})$ is compact, it is simple to show that, as $\varepsilon \to 0$,
$$II\leq \sup \left \{ \norm{\omega(e^{\alpha}\cdot)- \omega}_{H^s_r(\R^N)} \mid |\alpha| \leq \rho, \; \omega \in P_2(K_b^{PSP}) \right \}.$$ 
Summing up, we have
\begin{eqnarray*}
\dist((\lambda, u), K_b^{PSP})&\leq& e^{\frac{N\rho}{2}} \rho + \sup \left \{ \norm{\omega(e^{\alpha}\cdot)- \omega}_{H^s_r(\R^N)} \mid \abs{\alpha} \leq \rho, \; \omega \in P_2(K_b^{PSP}) \right \} \\
&\equiv& R(\rho)<\infty.
\end{eqnarray*}
Here we have
$$\lim_{\rho\to 0} R(\rho) =0,$$
which concludes the proof.
\QED

\bigskip

We are now ready to show that $\eta$ satisfies the desired properties.

\medskip

\claim Proof of Theorem \ref{thm_def_gen}.
Let $\mc{O}$ be a neighbourhood of $K_b^{PSP}$, and choose $R$ such that $N_R(K_b^{PSP})\subset \mc{O}$. 
By Lemma \ref{lem_intorni} choose $\rho\ll 1$ such that $R(\rho)<R$ and thus $N_{R(\rho)}(K_b^{PSP}) \subset \mc{O}$. 
Consequently, by Theorem \ref{thm_def_H}, there exists a deformation $\tilde{\eta}$ corresponding to the neighbourhood $\tilde{\mc{O}}:=\tilde{N}_{\rho}(\tilde{K}_b)$. 
We thus define $\eta$ by \eqref{eq_def_flow} and prove the properties. 
Start observing that
$$(\lambda, u)\in \mc{I}^{b\pm \delta} \implies b\pm \delta > \mc{I}(\lambda, u) = \mc{H}(\iota(\lambda, u)) \implies \iota(\lambda, u)\in \mc{H}^{b\pm \delta},$$
i.e. $\iota(\mc{I}^{b\pm \delta}) \subset \mc{H}^{b\pm \delta}$; similarly, $\pi(\mc{H}^{b\pm \delta})\subset \mc{I}^{b\pm \delta}$.
\begin{enumerate}
\item $\eta(0,\lambda, u) = \pi(\tilde{\eta}(0,\iota(\lambda,u))) = \pi(\iota(\lambda,u)) = (\lambda, u)$.
\item If $(\lambda,u)\in \mc{I}^{b-\bar\varepsilon}$, then $\iota(\lambda,u)\in\mc{H}^{b-\bar\varepsilon}$.
Thus $\eta(t,\lambda, u) = \pi(\tilde{\eta}(t,\iota(\lambda,u))) = \pi(\iota(\lambda, u))= (\lambda, u)$.
\item $\mc{I}(\eta(t,\lambda, u)) = \mc{I}(\pi(\tilde{\eta}(t,\iota(\lambda,u)))) = \mc{H}(\tilde{\eta}(t,\iota(\lambda,u))) \leq \mc{H}(\iota(\lambda,u)) = \mc{I}(\lambda, u)$.
\item If $K_b^{PSP}=\emptyset$, then $\tilde{K}_b=\emptyset$. 
Thus for $(\lambda,u)\in \mc{I}^{b+\varepsilon}$,
$\mc{I}(\eta(1,\lambda, u)) = \mc{I}(\pi(\tilde{\eta}(1,\iota(\lambda,u)))) = \mc{H}(\tilde{\eta}(1,\iota(\lambda,u))) \leq b-\varepsilon$.
\item We have, by previous arguments and Lemma \ref{lem_intorni}, that $\iota(\mc{I}^{b+\varepsilon}\setminus \mc{O}) 
= \iota(\mc{I}^{b+\varepsilon}\cap \complement \mc{O}) 
\subseteq \iota(\mc{I}^{b+\varepsilon}) \cap \iota(\complement \mc{O}) 
\subseteq \mc{H}^{b+\varepsilon} \cap (\complement \tilde{\mc{O}})
= \mc{H}^{b+\varepsilon}\setminus \tilde{\mc{O}}$ and thus
 $$\eta(1,\mc{I}^{b+\varepsilon}\setminus \mc{O})
 = \pi(\tilde{\eta}(1,\iota(\mc{I}^{b+\varepsilon}\setminus \mc{O}))) 
 \subset \pi(\tilde{\eta}(1,\mc{H}^{b+\varepsilon}\setminus \tilde{\mc{O}})) 
 \subset\pi(\mc{H}^{b-\varepsilon}) \subset \mc{I}^{b-\varepsilon}.$$ 
The other inclusion is similar and easier.
\item 
We write $\tilde \eta(t,\theta,\lambda,u)=\big(\tilde \eta_0(t,\theta,\lambda,u), \tilde \eta_1(t,\theta,\lambda,u), \tilde \eta_2(t,\theta,\lambda,u)\big)$. Then by definition
$$\big(\eta_1(t,\lambda, u), \eta_2(t,\lambda, u)\big) = \Big(\tilde{\eta}_1(t,0,\lambda, u), \tilde{\eta}_2\big(t,0,\lambda, u(e^{-\tilde{\eta}_0(t,0,\lambda,u)}\cdot)\big)\Big)$$
thus by the property 6 of Theorem \ref{thm_def_H},
\begin{eqnarray*}
\big(\eta_1(t,\lambda, -u), \eta_2(t,\lambda, -u)\big) &=& \Big(\tilde{\eta}_1(t,0,\lambda, -u), \tilde{\eta}_2\big(t,0,\lambda, -u(e^{-\tilde{\eta}_0(t,0,\lambda,-u)}\cdot)\big)\Big) \\
&=& \Big(\tilde{\eta}_1(t,0,\lambda, u), -\tilde{\eta}_2\big(t,0,\lambda, u(e^{-\tilde{\eta}_0(t,0,\lambda,u)}\cdot)\big)\Big) \\
&=&\big(\eta_1(t,\lambda, u), -\eta_2(t,\lambda, u)\big).
\end{eqnarray*}
\end{enumerate}
The theorem is hence proved.
\QED

\bigskip

Now we are ready to prove the main theorem for $\mc{H}$.

\medskip

\claim Proof of Theorem \ref{thm_def_H}.
To avoid cumbersome notation, we write $\xi=(\theta,\lambda, u)\in M$. Set 
$$M'=\{D\mc{H}(\xi)\neq 0\}.$$
It is known (see \cite{AM}) that there exists a pseudo-gradient on the Hilbert manifold $M$ associated to $\mc{H}$, namely a locally Lipschitz vector field $\mc{V}:M'\to TM$ such that
\begin{itemize}
\item[(a)] $\norm{\mc{V}(\xi)}_{\xi} \leq 2 \norm{D\mc{H}(\xi)}_{\xi,*}$;
\item[(b)] $D\mc{H}(\xi)\cdot \mc{V}(\xi) \geq \norm{D\mc{H}(\xi)}_{\xi,*}^2$;
\end{itemize}
in particular,
\begin{equation}\label{eq_pseudograd}
\frac{1}{2} \norm{\mc{V}(\xi)}_{\xi} \leq \norm{D\mc{H}(\xi)}_{\xi,*} \leq \norm{\mc{V}(\xi)}_{\xi}.
\end{equation}
Moreover, we can ask, in the construction of the pseudo-gradient, that $\mc{V}$ is $\G$-equivariant, since $\mc{H}$ is $\G$-invariant.
Namely, set $\mc{V}=(\mc{V}_0,\mc{V}_1,\mc{V}_2)$, then $\mc{V}_0$ and $\mc{V}_1$ are even in $u$, while $\mc{V}_2$ is odd in $u$.

By Corollary \ref{cor_PStilde}, there exists $\delta=\delta_{\frac{\rho}{3}}>0$ such that
	\begin{equation}\label{eq_cor_PSP_delta}
	\forall \, \xi \in \mc{H}^{b+\delta}_{b-\delta} \;\; s.t. \;\; \dist_M(\xi, \tilde{K}_b) > \frac{\rho}{3} \; : \; \norm{D\mc{H}(\xi)}_{\xi,*}> \delta .
	\end{equation}
We assume 
 \begin{equation}\label{choice-epsilon}
 \varepsilon < \min\Big\{ \frac{1}{2}\bar{\varepsilon}, \frac{1}{4} \delta, \frac16\rho\delta\Big\}.
 \end{equation}
Set the following
$$A=\mc{H}^{b+\varepsilon}_{b-\varepsilon}, \quad B=\mc{H}^{b+2\varepsilon}_{b-2\varepsilon}$$
and choose a locally Lipschitz function $g\in C(M,[0,1])$ such that
$$g=1 \; \hbox{ on $A$}, \quad g=0 \; \hbox{ on $\complement B$},$$
for instance $g(\xi)=\frac{d(\xi,\complement B)}{d(\xi,\complement B)+d(\xi,A)}$. \\
When $\tilde K_b\not=\emptyset$, we choose a locally Lipschitz function $\tilde g\in C(M,[0,1])$ such that
$$ \tilde g=0 \ \hbox{ on }\, \tilde N_{\frac{\rho}{3}}(\tilde K_b), \quad	\tilde g=1 \ \hbox{ on }\, \complement \tilde N_{\frac{2}{3}\rho}(\tilde K_b). $$
When $\tilde K_b=\emptyset$, we set $\tilde g\equiv 1$. 
Moreover we introduce
$$b(r)=\parag{\frac{1}{r} && \hbox{ if $r\geq 1$} \\ 1 && \hbox{ if $0\leq r < 1$}.}$$
Finally define
$$W(\xi)=-g(\xi) \tilde g(\xi)b\left(\norm{\mc{V}(\xi)}_{\xi}\right)\mc{V}(\xi)$$
and, fixed $\xi \in M$, consider the Cauchy problem
$$\parag{ \tilde{\eta}'& =& W(\tilde{\eta}), \\ \tilde{\eta}(0)&=&\xi.}$$
We have that $W$ is well defined on $M$ and
$$\norm{W(\xi)}_{\xi}\leq \norm{\mc{V}(\xi)}_{\xi} b\left(\norm{\mc{V}(\xi)}_{\xi}\right) \leq 1, $$
where we have used that $|g|$, $|\tilde g|$ $\leq 1$.
Therefore we have the global existence of a flow $\tilde{\eta}(t,\xi)$; we are interested in $\tilde{\eta}$ restricted to $[0,1]$.
We now verify the desired properties.
\begin{itemize}
\item[1)] $\tilde{\eta}(0,\xi)=\xi$ by construction of the flow.
\item[2)] If $\xi \in \mc{H}^{b-\bar{\varepsilon}}$, then $g(\xi)=0$, and thus $W(\xi)=0$. 
This means that $\tilde{\eta}(t,\xi)\equiv \xi$ is an equilibrium solution. 
Since $W\in Lip_{loc}(M)$ we have uniqueness of the solution, hence actually $\tilde{\eta}(t,\xi)\equiv \xi$.
\item[3)] We have
\begin{eqnarray*}
\lefteqn{\frac{d}{dt} \mc{H}(\tilde{\eta}(t,\xi)) 
= D\mc{H}(\tilde{\eta}(t,\xi)) \tilde{\eta}'(t,\xi)} \\
&=& -D\mc{H}(\tilde{\eta}(t,\xi)) \mc{V}(\tilde{\eta}(t,\xi)) g(\tilde{\eta}(t,\xi)) \tilde g(\tilde{\eta}(t,\xi))b\left(\norm{\mc{V}(\tilde{\eta}(t,\xi))}_{\tilde{\eta}(t,\xi)}\right)\\
&\leq& -\norm{D\mc{H}(\tilde{\eta}(t,\xi))}_{\tilde{\eta}(t,\xi),*}^2 g(\tilde{\eta}(t,\xi)) \tilde g(\tilde{\eta}(t,\xi)) b\left(\norm{\mc{V}(\tilde{\eta}(t,\xi))}_{\tilde{\eta}(t,\xi)}\right) \\
&\leq& 0
\end{eqnarray*}
that is the claim; we have used that $g$, $\tilde g$, $b$ are positive and the property $(b)$.

\item[4)]
We assume here $\tilde{K}_b= \emptyset$.
By using the fundamental theorem of calculus and previous arguments, we obtain
\begin{eqnarray*}
\lefteqn{\mc{H}(\tilde{\eta}(1,\xi)) - \mc{H}(\tilde{\eta}(0,\xi)) = \int_0^1 \frac{d}{ds} \mc{H}(\tilde{\eta}(s,\xi)) ds} \\
&=& -\int_0^1 D\mc{H}(\tilde{\eta}(s,\xi)) \mc{V}(\tilde{\eta}(s,\xi)) g(\tilde{\eta}(s,\xi)) b\left(\norm{\mc{V}(\tilde{\eta}(s,\xi))}_{\tilde{\eta}(s,\xi)}\right) ds \\
&\leq& -\int_0^1 \norm{D\mc{H}(\tilde{\eta}(s,\xi))}_{\tilde{\eta}(s,\xi),*}^2 g(\tilde{\eta}(s,\xi)) b\left(\norm{\mc{V}(\tilde{\eta}(s,\xi))}_{\tilde{\eta}(s,\xi)}\right) ds.
\end{eqnarray*}
Let now $\xi \in \mc{H}^{b+\varepsilon}$. This means, by point $3)$, that for $s\in [0,1]$
$$\mc{H}(\tilde{\eta}(s,\xi))\leq \mc{H}(\tilde{\eta}(0,\xi)) = \mc{H}(\xi) \leq b + \varepsilon,$$
thus $\tilde{\eta}(s,\xi)\in \mc{H}^{b+\varepsilon}$ and
$$ \mc{H}(\tilde{\eta}(1,\xi)) \leq b + \varepsilon - \int_0^1 \norm{D\mc{H}(\tilde{\eta}(s,\xi))}_{\tilde{\eta}(s,\xi),*}^2 g(\tilde{\eta}(s,\xi)) b\left(\norm{\mc{V}(\tilde{\eta}(s,\xi))}_{\tilde{\eta}(s,\xi)}\right) ds. $$
Assume now by contradiction that $\mc{H}(\tilde{\eta}(1,\xi))> b- \varepsilon$, which implies (again by point $3)$) $\mc{H}(\tilde{\eta}(s,\xi))> b- \varepsilon$, for all $s\in [0,1]$. 
Thus for all $s\in [0,1]$ we have $\tilde{\eta}(s,\xi) \in \mc{H}^{b+\varepsilon}_{b-\varepsilon}$ and in particular, since $\varepsilon < \frac{1}{2} \bar{\varepsilon}$, that $g(\tilde{\eta}(s,\xi))=1$; hence
$$ \mc{H}(\tilde{\eta}(1,\xi)) \leq b + \varepsilon - \int_0^1 \norm{D\mc{H}(\tilde{\eta}(s,\xi))}_{\tilde{\eta}(s,\xi),*}^2 b\left(\norm{\mc{V}(\tilde{\eta}(s,\xi))}_{\tilde{\eta}(s,\xi)}\right) ds. $$
By \eqref{eq_pseudograd}, by the fact that $\tilde{\eta}(s,\xi) \in \mc{H}^{b+\varepsilon}_{b-\varepsilon}\subset \mc{H}^{b+\delta}_{b-\delta}$ and by \eqref{eq_curv_pend_gen}, we have
\begin{equation}\label{eq_dim_stimaV}
\norm{\mc{V}(\tilde{\eta}(s,\xi))}_{\tilde{\eta}(s,\xi)} \geq \norm{D\mc{H}(\tilde{\eta}(s,\xi))}_{\tilde{\eta}(s,\xi),*} \geq \delta \geq 4 \varepsilon;
\end{equation}
in particular, 
$$ b\left(\norm{\mc{V}(\tilde{\eta}(s,\xi))}_{\tilde{\eta}(s,\xi)}\right) =\frac{1}{\norm{\mc{V}(\tilde{\eta}(s,\xi))}_{\tilde{\eta}(s,\xi)}}.$$
Thus, exploiting again \eqref{eq_pseudograd} and \eqref{eq_dim_stimaV} we have
\begin{eqnarray*}
\mc{H}(\tilde{\eta}(1,\xi)) &\leq& b + \varepsilon - \frac{1}{2} \int_0^1 \norm{D\mc{H}(\tilde{\eta}(s,\xi))}_{\tilde{\eta}(s,\xi),*} ds\\
 &\leq& b + \varepsilon - 2\int_0^1 \varepsilon ds = b - \varepsilon,
\end{eqnarray*}
which is an absurd.
\item[5)] We assume now $\tilde{K}_b \neq \emptyset$.
Let now $\xi \in \mc{H}^{b+\varepsilon}\setminus \tilde{\mc{O}}$. 
Assume again by contradiction that $\mc{H}(\tilde{\eta}(1,\xi))>b-\varepsilon$, which implies again $\tilde{\eta}(s,\xi)\in \mc{H}^{b+\varepsilon}_{b-\varepsilon}$. We distinguish two cases.
\\ \textbf{Case 1:} $\tilde{\eta}(t,\xi)\notin \tilde{N}_{\frac{2}{3}\rho}(\tilde{K}_b)$ for all $t\in [0,1]$. 
In this case we proceed as in the proof of 4). Indeed since $\varepsilon<\delta_{\frac{\rho}{3}}$, we are in the case of \eqref{eq_curv_pend_gen} and thus
$$\norm{D\mc{H}(\tilde{\eta}(s,\xi))}_{\tilde{\eta}(s,\xi),*}> \delta>4 \varepsilon.$$
We can do exactly the same passages as before and conclude.
\\ \textbf{Case 2:} $\tilde{\eta}(t^*,\xi)\in \tilde{N}_{\frac{2}{3}\rho}(\tilde{K}_b)$ for some $t^*\in [0,1]$. 
In this case $\varepsilon$ has to be better specified. We make a finer argument by choosing suitable $[\alpha,\beta]\subset [0,1]$ and observing that
\begin{eqnarray*}
 \mc{H}(\tilde{\eta}(1,\xi)) &\leq& \mc{H}(\tilde{\eta}(\beta,\xi)) = \mc{H}(\tilde{\eta}(\alpha,\xi)) + \int_{\alpha}^{\beta} \frac{d}{ds} \mc{H}(\tilde{\eta}(s,\xi)) ds \\
&\leq& \mc{H}(\tilde{\eta}(0,\xi)) + \int_{\alpha}^{\beta} \frac{d}{ds} \mc{H}(\tilde{\eta}(s,\xi)) ds \\
&\leq& b+\varepsilon + \int_{\alpha}^{\beta} \frac{d}{ds} \mc{H}(\tilde{\eta}(s,\xi)) ds.
\end{eqnarray*}
Noting that $\tilde{\eta}(0,\xi)=\xi\notin \tilde{\mc{O}}=\tilde{N}_{\rho}(\tilde{K}_b)$ and $\tilde{\eta}(t^*,\xi)\in \tilde{N}_{\frac{2}{3}\rho}(\tilde{K}_b)$, we can find $\alpha$ and $\beta$ such that
$$\tilde{\eta}(\alpha)\in \partial \tilde{N}_{\rho}(\tilde{K}_b), \quad \tilde{\eta}(\beta)\in \partial \tilde{N}_{\frac{2}{3}\rho}(\tilde{K}_b),$$
and
$$\tilde{\eta}(s)\in \tilde{N}_{\rho}(\tilde{K}_b) \setminus \tilde{N}_{\frac{2}{3}\rho}(\tilde{K}_b) \quad \forall \, s\in (\alpha,\beta).$$
Hence we obtain by \eqref{eq_cor_PSP_delta}
$$\mc{H}(\tilde{\eta}(1,\xi)) \leq b+\varepsilon - \delta(\beta-\alpha).$$
We need an estimate from below of $\beta-\alpha$. We easily obtain it by observing that $\tilde{\eta}(\cdot,\xi)$ is a path connecting $\tilde{\eta}(\alpha, \xi)$ and $\tilde{\eta}(\beta, \xi)$, thus (recall that $1\geq \norm{W(\xi)}_{\xi}$)
\begin{eqnarray*}
\beta - \alpha &=& \int_{\alpha}^{\beta} dt \geq \int_{\alpha}^{\beta} \norm{W(\tilde{\eta}(t,\xi))}_{\tilde{\eta}(t,\xi)} dt \\
&=& \int_{\alpha}^{\beta} \norm{\tilde{\eta}'(t,\xi)}_{\tilde{\eta}(t,\xi)} dt \geq \dist_M(\tilde{\eta}(\alpha, \xi),\tilde{\eta}(\beta,\xi)) \\
&\geq&\dist_M\left(\tilde{N}_{\rho}(\tilde{K}_b), \tilde{N}_{\frac{2}{3}\rho}(\tilde{K}_b)\right) \geq \frac{1}{3}\rho.
\end{eqnarray*}
Finally
$$\mc{H}(\tilde{\eta}(1,\xi)) \leq b+\varepsilon - \tfrac{1}{3} \rho \delta\leq b-\varepsilon$$
by our choice \eqref{choice-epsilon} of $\varepsilon$.

As regards the second inclusion, we argue in a similar way. Let $\xi \in \mc{H}^{b+\varepsilon}$. Case 1 can be done verbatim. 
In Case 2, if $\tilde{\eta}(1,\xi)\in \tilde{\mc{O}}$ we are done; if not, then we repeat the argument but with the path built thanks to $\tilde{\eta}(1,\xi)\notin \tilde{N}_{\rho}(\tilde{K}_b)$ and $\tilde{\eta}(t^*,\xi)\in \tilde{N}_{\frac{2}{3}\rho}(\tilde{K}_b)$.
\item[6)] Notice that, written $W=(W_0,W_1,W_2)$, we have that $W_0$ and $W_1$ are even in $u$ while $W_2$ is odd in $u$, since $\mc{V}$ is so and $g$, $b(\norm{D\mc{H}(\cdot)}_{\cdot,*})$ are even in $u$. 
Thus, by uniqueness of the solution, we have that $\tilde{\eta}$ satisfies the required symmetry properties.
\end{itemize}
The proof is thus concluded.
\QED


\subsection{Minimax values $a_j(\lambda)$}

We write for $j\in\N$, $D_j=\{\xi\in\R^j \mid \abs\xi\leq 1\}$ and we introduce the set of paths
$$\Gamma_j(\lambda)=\big\{\gamma\in C(D_j, H^s_r(\R^N)) \, \mid \, \gamma \hbox{ odd}, \, \mc{J}(\lambda, \gamma(\xi))<0 \; \forall \xi \in \partial D_j \big\}$$
and
$$a_j(\lambda)=\inf_{\gamma \in \Gamma_j(\lambda)} \sup_{\xi \in D_j} \mc{J}(\lambda, \gamma(\xi)).$$
By an odd extension from $[0,1]$ to $[-1,1]=D_1$, we may regard $\Gamma_1(\lambda)\equiv\Gamma(\lambda)$ and $a_1(\lambda)\equiv a(\lambda)$. 
Thus these quantities are a good generalization. As for $j=1$, we prove the following properties.

\begin{proposition}
Let $\lambda_0\in (-\infty,\infty]$ be the number given in \eqref{muzero}--\eqref{lambdazero}, $\lambda < \lambda_0$ and $j \in \N$.
\begin{enumerate}
\item $\Gamma_j(\lambda)\neq \emptyset$, thus $a_j(\lambda)$ is well defined. Moreover, it is increasing with respect to $\lambda$;
\item $a_j(\lambda)\leq a_{j+1}(\lambda)$;
\item $a_j(\lambda)> 0$; 
\item $\lim_{\lambda\to \lambda_0^-} \frac{a_j(\lambda)}{e^{\lambda}}=+\infty$;
\item if \textnormal{(g4)} holds, then $\lim_{\lambda \to -\infty} \frac{a_j(\lambda)}{e^{\lambda}}=0$.
\end{enumerate}
\end{proposition}


\claim Proof.
The proofs are quite the same of Propositions \ref{lem_buona_def}--\ref{out}. We point out just some slight differences.
\begin{enumerate}
\item For $\lambda<\lambda_0$, there exists $\xi_0>0$ such that 
	$$ G(t_0)-{e^\lambda\over 2}t_0^2 >0. $$
As in \cite{BL2}, we find that there exists a continuous odd map $\hat{\gamma}:\,\partial D_j\to H_r^1(\R^N)\hookrightarrow H^s_r(\R^N)$ with $\mc{J}(\lambda,\hat{\gamma}(\xi))<0$. 
Extending $\hat{\gamma}$ onto $D_j$ we find $\Gamma_j(\lambda)\not=\emptyset$.
\item Since $D_j\subset D_{j+1}$, we observe $\gamma_{|D_j}\in\Gamma_j(\lambda)$ 
for $\gamma\in\Gamma_{j+1}(\lambda)$. Thus we regard $\Gamma_{j+1}(\lambda)\subset\Gamma_j(\lambda)$ 
and obtain 2).
\item Clear by $a_1(\lambda)=a(\lambda)>0$ and point $2)$. 
\item Again by $\lim_{\lambda\to \lambda_0^-} \frac{a(\lambda)}{e^{\lambda}}=+\infty$ and point $2)$.
\item We consider the path $\hat{\gamma}:\partial D_j\to H^s_r(\R^N)$ obtained in 1) and introduce a path
$$ \xi \mapsto \mu^{N/4} \hat{\gamma}\left({\xi\over\abs\xi}\right)(\cdot/\mu^{-1/(2s)}\abs\xi); D_j \to H^s_r(\R^N).
$$
Arguing as in Proposition \ref{out}, we have 5).
\QED
\end{enumerate}


\subsection{Minimax values $b^m_j$}

We set
\begin{eqnarray*}
\Gamma_j^m=\{\Theta\in C(D_j, \R\times H^s_r(\R^N)) &\mid&
 \hbox{$\Theta$ is $\G$-equivariant;} \\
&& \mc{I}(\Theta(0)) \leq B_m-1; \\
&& \Theta(\xi)\notin \Omega, \ \mc{I}(\Theta(\xi))\leq B_m-1 \ \hbox{for all}\ \xi \in \partial D_j\}
\end{eqnarray*}
and
$$b_j^m= \inf_{\Theta \in \Gamma_j^m} \sup_{\xi \in D_j} \mc{I}(\Theta(\xi)).$$
We notice that for $j=1$ we obtain $\Gamma_1^m \equiv \Gamma^m$ (up to an even/odd extension from $[0,1]$ to $[-1,1]=D_1$) and $b_1^m \equiv b_m$. 
So $\Gamma_j^m$ is the natural extension to build multiple solutions.

As in the case of $\Gamma^m$ and $b_m$, we want to prove that $\Gamma_j^m \neq \emptyset$ and that, for a fixed $k\in \N$, there exists an $m_k \gg 0$ (possibly equal to $0$) such that, if $m>m_k$, then $ b_j^m <0$ for $j=1 \dots k$.


\begin{proposition}\label{prop_mk}
$\,$
\begin{itemize}
\item[\textnormal{(i)}] For any $\lambda <\lambda_0$, $m>0$, $j \in \N$, we have $\Gamma_j^m\neq \emptyset$ and $b_j^m\leq a_j(\lambda) - e^{\lambda} \frac{m}{2}$.
\item[\textnormal{(ii)}] For any $k\in\N$ there exists $m_k\geq 0$ such that for $m>m_k$
	$$ b_j^m < 0 \quad \hbox{for}\ j=1,2,\dots, k. $$
\item[\textnormal{(iii)}] $m_k=0$ for all $k\in\N$ if \textnormal{(g4)} holds. That is,
	$$ b_j^m < 0 \quad \hbox{for all}\ j\in\N. $$
\end{itemize}
\end{proposition}

\claim Proof. For (i), the proof is similar to Proposition \ref{tom}. We just need to set
$\zeta_{\lambda}\in \Gamma_j(\lambda)$, 
	$$ \psi_{\lambda}(\xi) = 
	\parag{
		(\lambda+L(1-2|\xi|), \,0) && \hbox{ if $|\xi| \in [0,1/2]$,} \\ 
		\left(\lambda, \, \zeta_{\lambda}\left({\xi\over\abs{\xi}}(2\abs\xi-1)\right)\right) && \hbox{ if $|\xi| \in (1/2,1]$}
		} $$	
and we come up again to the same proof.

For (ii), (iii), with a proof similar to Proposition \ref{blue}, we set 
\begin{equation}\label{defmk}
	m_k= 2\inf_{\lambda<\lambda_0}\frac{a_k(\lambda)}{e^{\lambda}}\geq 0
\end{equation}
and observe that $m_k\leq m_{k+1}$ since $a_k(\lambda)$ are increasing in $k$.
\QED


\subsection{Minimax values $c_j^m$}

Let us define minimax families $\Lambda_j^m$ which allow to find multiple solutions. We use an idea from \cite{Rab0}.
In what follows, we denote by $\genus(A)$ the genus of closed symmetric sets $A$ with $0\not\in A$. 

Define, for each $j\in \N$, 
\begin{eqnarray*}
\Lambda_j^m=\{A=\Theta(\overline{D_{j+l}\setminus Y}) &\mid & l\geq 0, \; \Theta\in \Gamma_{j+l}^m, \\
&& Y\subseteq D_{j+l}\setminus \{0\} \; \hbox{ is closed, symmetric in $0$ and $\genus(Y)\leq l$} \}
\end{eqnarray*}
and
$$c_j^m= \inf_{A\in \Lambda_j^m} \sup_{A} \mc{I}.$$

In the following lemma, we observe that $\Lambda_j^m$ includes, in some way, $\Gamma^m_j$ and that it inherits the property that the paths intersect $\partial \Omega$.

\begin{lemma}\label{lem7.6}
$\,$
\begin{itemize}
\item[\textnormal{(i)}] $\Lambda_j^m \neq \emptyset$;
\item[\textnormal{(ii)}] $c_j^m\leq b_j^m$;
\item[\textnormal{(iii)}]
for any $A\in \Lambda_j^m$, we have $A\cap \partial \Omega \neq \emptyset$. As a consequence, we obtain
$$b_m=B_m= B'_m\leq c_j^m.$$
\end{itemize}
\end{lemma}

\claim Proof.
Indeed, we see that, by choosing $l=0$ and $Y=\emptyset$ we have
$$\{A=\Theta(D_j) \, \mid \, \Theta \in \Gamma_j^m\} \subset \Lambda_j^m$$
from which easily come the first two claims.

Focus on the third claim. Let $A=\Theta(\overline{D_{j+l}\setminus Y})$ and set $U=\Theta^{-1}(\Omega)$. 
By the symmetry in $(\lambda, u)$ of $\Theta$ and the symmetry in $u$ of $\Omega$ we have that $U$ is symmetric. 
Moreover, since $\Theta(0)\in \Omega$, we have that $U\subset D_{j+l} \subset \R^{j+l}$ is a symmetric neighbourhood of the origin. 
By the properties of the genus we have that the genus of $\partial U$ is maximum, that is
\begin{equation}\label{eq_genus}
\genus(\partial U) = j+l.
\end{equation}
Observe in addition the following chain of inclusions
$$\overline{\partial U \setminus Y} = \overline{(\partial U \cap D_{j+l})\setminus Y} = \overline{(D_{j+l} \setminus Y) \cap \partial U} \subseteq \overline{D_{j+l} \setminus Y} \cap \overline{\partial U} = \overline{D_{j+l}\setminus Y} \cap \partial U $$
thus
$$\Theta\left(\overline{\partial U \setminus Y}\right) \subseteq \Theta\left( \overline{D_{j+l}\setminus Y} \cap \partial U\right) \subseteq \Theta\left( \overline{D_{j+l}\setminus Y}\right) \cap \Theta \left(\partial U\right) = A\cap \Theta \left(\partial U\right)$$
and since $\partial U = \partial \left( \Theta^{-1}(\Omega)\right) \subseteq \Theta^{-1}(\partial \Omega)$ we obtain
$$\Theta\left(\overline{\partial U \setminus Y}\right) \subseteq A \cap \partial \Omega. $$
Thus, to show the claim, it is sufficient to show that $\overline{\partial U \setminus Y}\neq \emptyset$. 
But is an immediate consequence of \eqref{eq_genus} and the property of the genus
$$ \genus(\overline{\partial U \setminus Y}) \geq \genus(\partial U) - \genus(Y) \geq (j+l)-l = j \geq 1$$
which directly excludes the possibility that $\overline{\partial U \setminus Y}$ is empty. This concludes the proof of the first part.

We prove now the consequence. Indeed, for each $A\in \Lambda_j^m$ we have
$$B_m' = \inf_{\partial \Omega} \mc{I} \leq \inf_{\partial \Omega \cap A} \mc{I} \leq \sup_{\partial \Omega \cap A} \mc{I} \leq \sup_{A} \mc{I}$$
and thus the claim passing to the infimum over $\Lambda_j^m$.
\QED

\bigskip

Let us now show the main properties of $\Lambda_j^m$ and $c_j^m$. 
We point out that these classical properties are the only ones which will be used in the proof of the existence of multiple solutions.

\begin{proposition} \label{prop_genus}$\;$
\begin{enumerate}
\item $\Lambda_j^m \neq \emptyset$;
\item $\Lambda_{j+1}^m\subseteq \Lambda_j^m$, and thus $c_j^m\leq c_{j+1}^m$;
\item let $A\in \Lambda_j^m$ and $Z\subset \R \times H^s_r(\R^N)$ be $\G$-invariant, closed, 
and such that $0 \notin \overline{P_2(Z)}$ and $\genus(\overline{P_2(Z)})\leq i$. Then $\overline{A\setminus Z} \in \Lambda_{j-i}^m$.
\end{enumerate}
Fix now $k\in \N$, and let $m>m_k$, where $m_k$ has been introduced in Proposition \ref{prop_mk}, i.e., in \eqref{defmk}.
Then
\begin{itemize}
\item[4.] $c_j^m<0$ and $\mc{I}$ satisfies $(PSP)_{c_j^m}$;
\item[5.] if $A\in \Lambda_j^m$ and $\eta$ is a deformation as in Theorem \ref{thm_def_gen} for $b=c_j^m$, then $\eta(1,A)\in \Lambda_j^m$.
\end{itemize}
\end{proposition}


\claim Proof.
Properties $1)$ and $4)$ has already been shown in the Lemma \ref{lem7.6}, while property $2)$ is a consequence of the definition. 
Let us see properties $3)$ and $5)$.
\begin{itemize}
\item[3)]
Let $A=\Theta(\overline{D_{j+l}\setminus Y}) \in \Lambda_j^m$ and let $Z$ be $\G$-invariant, closed and such that $0\notin \overline{P_2(Z)}$ and $\genus(\overline{P_2(Z)})\leq i$. 
Assume it holds
\begin{eqnarray}\label{eq_Gamma_g}
\overline{A\setminus Z} &=& \Theta((\overline{D_{j+l}\setminus Y)\setminus \Theta^{-1}(Z)}) \\
&=& \Theta(\overline{D_{(j-i)+(l+i)} \setminus (Y \cup \Theta^{-1}(Z))}); \notag
\end{eqnarray}
if $\genus(Y \cup \Theta^{-1}(Z))\leq l+i$ we have the claim. But this is a direct consequence of the assumptions and the property of the genus; indeed
\begin{eqnarray*}
\genus(Y \cup \Theta^{-1}(Z))&\leq& \genus(Y)+\genus(\Theta^{-1}(Z)) \\
&\leq& l + \genus(\overline{h(\Theta^{-1}(Z))}) \\
&=& l + \genus(\overline{P_2(Z)}) \\
&\leq& l+i
\end{eqnarray*}
where we have set $h=P_2 \circ \Theta$, which is an odd map and thus admissible for the genus.

Turn now to \eqref{eq_Gamma_g}. Set $B=D_{j+l}\setminus Y$ and $W=\Theta^{-1}(Z)$ we have to prove
$$\overline{\Theta(\overline{B})\setminus \Theta(W)} = \Theta(\overline{B\setminus W}).$$
We have
$$\overline{\Theta(\overline{B})\setminus \Theta(W)} \subseteq
 \overline{\Theta(\overline{B}\setminus W)} \stackrel{(i)} \subseteq \overline{\Theta(\overline{B\setminus W})} \stackrel{(ii)}= \Theta(\overline{B\setminus W}) $$
and
$$\Theta(\overline{B\setminus W}) \stackrel{(iii)} \subseteq \overline{\Theta(B\setminus W)} \stackrel{(iv)} = \overline{\Theta(B)\setminus \Theta(W)} \subseteq \overline{\Theta(\overline{B})\setminus \Theta(W) }$$
where
\begin{itemize}
\item[\textnormal{(i)}] is due to the fact that $W$ is closed;
\item[\textnormal{(ii)}] $\overline{B\setminus W}\subseteq D_{j+l}$ is compact, thus $\Theta(\overline{B\setminus W})$ is closed;
\item[\textnormal{(iii)}] derives from the continuity of $\Theta$;
\item[\textnormal{(iv)}] is due to the fact that $W$ is a preimage.
\end{itemize}
\item[5)]
Consider $0<\bar{\varepsilon}<1$, $b=c_j^m\geq B_m$ and $\eta$ as in the deformation lemma, and fix $A=\Theta(\overline{D_{j+l}\setminus Y}) \in \Lambda_j^{m}$ with $\Theta \in \Gamma_{j+l}^m$. 
To show that $\eta(1,A)\in \Lambda_j^m$ and conclude the proof, it is sufficient to show that $\tilde{\Theta}=\eta(1, \Theta)\in \Gamma_{j+l}^m$ as well.
\begin{itemize}
\item[$\bullet$] $\tilde{\Theta}(-\xi) = \eta(1,\Theta(-\xi)) = \eta(1,\Theta_1(-\xi), \Theta_2(-\xi)) = \eta(1,\Theta_1(\xi), -\Theta_2(\xi))$
and thus
\begin{eqnarray*}
\lefteqn{\big(\tilde{\Theta}_1(-\xi), \tilde{\Theta}_2(-\xi)\big) = \Big(\eta_1\big(1,\Theta_1(\xi), -\Theta_2(\xi)\big), \eta_2\big(1,\Theta_1(\xi), -\Theta_2(\xi)\big)\Big)} \\
&=& \Big(\eta_1\big(1,\Theta_1(\xi), \Theta_2(\xi)\big), -\eta_2\big(1,\Theta_1(\xi), \Theta_2(\xi)\big)\Big) 
= \big(\tilde{\Theta}_1(\xi), -\tilde{\Theta}_2(\xi)\big) 
\end{eqnarray*}
which shows that $\tilde{\Theta}_1$ is even and $\tilde{\Theta}_2$ is odd.
\item[$\bullet$] By Lemma \ref{lem7.6}, for $\xi=0$ and $\xi \in \partial D_{j+l}$ we have $\mc{I}(\Theta(\xi)) \leq B_m-1 = B'_m-1 \leq c_j^m-\bar{\varepsilon}$, thus $\Theta(\xi)\in \mc{I}^{c_j^m-\bar{\varepsilon}}$. 
Therefore $\tilde{\Theta}(\xi)=\eta(1,\Theta(\xi))=\Theta(\xi)$ for $\xi=0$ and $\xi \in \partial D_{j+l}$, and the same properties are satisfied.
\QED
\end{itemize}
\end{itemize}


\subsection{Multiplicity theorem}

Fix $k\in \N$, and let $\Lambda_j^m$ and $c_j^m$ be given in the previous section.
By the properties given in Proposition \ref{prop_genus}, we can find multiple solutions.


\begin{theorem}\label{theopassmult}
Fix $k\in \N$, and assume $m>m_k$. We have that
$$c_1^m \leq c_2^m \leq \dots \leq c_k^m<0$$
are critical values of $\mc{I}$. Moreover
\begin{itemize}
\item[\textnormal{(i)}]
 if, for some $q\geq 1$,
$$c_j^m < c_{j+1}^m < \dots < c_{j+q}^m$$
then we have $q+1$ different nonzero critical values, and thus $q+1$ different pairs of nontrival solutions of \eqref{problem};
\item[\textnormal{(ii)}]
if instead, for some $q\geq 1$,
\begin{equation}\label{eq7.6}
	c_j^m = c_{j+1}^m = \dots = c_{j+q}^m \equiv b
\end{equation}
then 
\begin{equation}\label{eq_genus_teo}
\genus(P_2(K_b^{PSP}))\geq q+1
\end{equation}
and thus, by the properties the genus, $\# P_2(K_b^{PSP})=+\infty$, which means that we have infinite different solutions of \eqref{problem}. %
\end{itemize}%
Summing up, we have at least $k$ different (pairs of) solutions of \eqref{problem}.
\end{theorem}

\claim Proof.
It is sufficient to show only the property \eqref{eq_genus_teo} on the genus: indeed by choosing $q=0$ we have that, for each $j$, $\#(K_{c_j^m}^{PSP})\geq 1$ and thus $c_j^m$ is a nontrivial critical value.

By the $(PSP)_b$ we have that $K_b^{PSP}$ is compact, thus $P_2(K_b^{PSP})$ is compact; moreover it is symmetric with respect to $0$ and does not contain $0$ (see Corollary \ref{PSP-cor}).

By the property of the genus we can find a (closed, symmetric with respect to origin, not containing the zero) neighbourhood $N$ of $P_2(K_b^{PSP})$ which preserves the genus, i.e. $\genus(N)=\genus(P_2(K_b^{PSP}))$. 
We can easily think $N$ as a projection of a neighbourhood $Z$ of $K_b^{PSP}$ (i.e. $N=P_2(Z)$) satisfying the properties of Proposition \ref{prop_genus}. 

By Theorem \ref{thm_def_gen}, there exist a sufficiently small $\varepsilon$ and an $\eta$ such that $\eta(\mc{I}^{b+\varepsilon}\setminus Z) \subseteq \mc{I}^{b-\varepsilon}$. 
Corresponding to $\varepsilon$, by definition of $c_j^m$, there exists an $A\in \Lambda_{j+q}^m$ such that $\sup_{A}\mc{I}<b+\varepsilon$, that is $A\subseteq \mc{I}^{b+\varepsilon}$. Thus 
$$ \eta(1, \overline{A\setminus Z}) \subseteq \eta(1, \overline{\mc{I}^{b+\varepsilon}\setminus Z}) \stackrel{\eta(1,\cdot) \hbox{ continuous}} \subseteq \overline{\eta(1,\mc{I}^{b+\varepsilon}\setminus Z)} \subseteq \overline{\mc{I}^{b-\varepsilon}} = \mc{I}^{b-\varepsilon},$$
and hence
\begin{equation}\label{eq_th_molt}
\sup_{\eta(1,\overline{A\setminus Z})} \mc{I} \leq b-\varepsilon.
\end{equation}
On the other hand, assume by contradiction that $\genus(P_2(K_b^{PSP}))\leq q$, i.e. $\genus(P_2(Z))\leq q$. We use now the properties on $c_j^m$ and $\Lambda_j^m$.

Replacing $j$ with $j+q$ and $i$ with $q$ and applying Proposition \ref{prop_genus}, we have $\overline{A\setminus Z}\in \Lambda_j^m$; by property 5) of Proposition \ref{prop_genus} we obtain $\eta(1,\overline{A\setminus Z}) \in \Lambda_j^m$, which implies (by definition of $c_j^m$) 
$$\sup_{\eta(1,\overline{A\setminus Z})} \mc{I} \geq c_j^m=b.$$
This is a contradiction with \eqref{eq_th_molt}, and thus concludes the proof.
\QED

\bigskip
\claim Proof of Theorem \ref{S:1.13}.
As consequence of Theorem \ref{theopassmult}, we derive (i). 
We pass to prove (ii). Under condition \textnormal{(g4)}, we have $m_k=0$ for all $k \in \N$. 
Thus for any $j \in \N$, $c_j^m$ is a critical value of $\mathcal{I}$ and $c_j^m \leq b_j^m <0$. Since $c_j^m$ is an increasing sequence, we have $c_j^m \to \bar c \leq 0$ as $j \to \infty$. We need to show that $\bar c=0$.

By contradiction we assume $\bar c <0$. Then $K^{PSP}_{\bar c}$ is compact and $K^{PSP}_{\bar c} \cap (\R \times \{0\}) =\emptyset$. It follows that $q= \genus(P_2(K^{PSP}_{\bar c})) < \infty$.
Arguing as in the proof of Theorem \ref{theopassmult}, let $\delta>0$ such that $q= \genus(P_2(N_\delta(K^{PSP}_{\bar c}))) < \infty$.
By Theorem \ref{thm_def_gen}, there exist $\varepsilon\in (0,1)$ small and $\eta: [0,1] \times \R \times H^s_r(\R^N) \to \R \times H^s_r(\R^N)$ satisfying
\begin{equation}\label{ax1}
\eta(1,\mc{I}^{\bar c+ \varepsilon}\setminus N_\delta(K^{PSP}_{\bar c})) \subseteq \mc{I}^{ \bar c-\varepsilon}
\end{equation}
and
\begin{equation}\label{ax2}
\eta(t,\lambda, u) =(\lambda,u) \quad \hbox{if $\, \mc{I}(\lambda, u) \leq B_m -1$}. 
\end{equation} 
We can choose $j \in \N$ sufficiently large such that $c_j^m > \bar c - \varepsilon$ and take $B \in \Lambda_{j+q}^m$ such that $B \subset \mc{I}^{\bar c +\varepsilon}$.
Then we have 
$$ \overline{B \setminus N_\delta(K^{PSP}_{\bar c})} \in \Lambda_j^m.$$
From equations \eqref{ax1}, \eqref{ax2} we derive $c_j^m \leq \bar c - \varepsilon$, which gives a contradiction.
\QED

\bigskip

\noindent {\bf Acknowledgments.} 
The first and second authors are supported by PRIN 2017JPCAPN ``Qualitative and quantitative aspects of nonlinear PDEs'', and partially supported by GNAMPA-INdAM.
The third author is supported in part by Grant-in-Aid for Scientific Research (JP19H00644, JP18KK0073, JP17H02855, JP16K13771, JP26247014) of Japan Society for the Promotion of Science.

%
%


\begin{thebibliography} {l}

	\bibitem{AM} A. Ambrosetti, A. Malchiodi, 
	``Nonlinear analysis and semilinear elliptic problems'', 
	Cambridge Univ. Press, New York, 2010.

	\bibitem{AR} A. Ambrosetti, P. H. Rabinowitz, 
	\emph{Dual variational methods in critical point theory and applications},
	J. Funct. Anal. \textbf{14} (1973), no. 4, 349--381.

	\bibitem{Ata0} A. Atangana, 
	``Fractional operators with constant and variable order with application to geo-hydrology'',
	Academic Press, London, 2018.
	
	\bibitem{BV} T. Bartsch, S. de Valeriola,
	\emph{Normalized solutions of nonlinear Schr\"odinger equations},
	Arch. Math. (Basel) \textbf{100} (2013), no. 1, 75--83.
	
	\bibitem{BGMMV} J. Bellazzini, M. Ghimenti, C. Mercuri, V. Moroz, J. Van Schaftingen,
	\emph{Sharp Gagliardo-Nirenberg inequalities in fractional Coulomb-Sobolev spaces},
	Trans. Amer. Math. Soc. \textbf{370} (2018), no. 11, 8285--8310.
	
	\bibitem{BL1} H. Berestycki, P.-L. Lions,
	\emph{Nonlinear scalar field equations. I. Existence of a ground state},
	Arch. Rational Mech. Anal. \textbf{82} (1983), no. 4, 313--345.
	
	\bibitem{BL2} H. Berestycki, P.-L. Lions,
	\emph{Nonlinear scalar field equations. II. Existence of infinitely many solutions},
	Arch. Rational Mech. Anal. \textbf{82} (1983), no. 4, 347--375.
	
	\bibitem{BKS} J. Byeon, O. Kwon, J. Seok,
	\emph{Nonlinear scalar field equations involving the fractional Laplacian},
	Nonlinearity \textbf{30} (2017), no. 4, 1659--1681.

 	\bibitem{BJM} J. Byeon, L. Jeanjean, M. Maris,
 	\emph{Symmetry and monotonicity of least energy solutions},
 	Calc. Var. Partial Differential Equations \textbf{36} (2009), no. 4, 481--492.
 		
	\bibitem{BuV} C. Bucur, E. Valdinoci,
	``Nonlocal Diffusion and Applications'', 
	Lect. Notes Unione Mat. Ital. \textbf{20}, Springer, Switzerland, 2016.

	\bibitem{CL} T. Cazenave, P.-L. Lions, 
	\emph{Orbital stability of standing waves for some nonlinear Schr\"odinger equations}, 
	Commun. Math. Phys. \textbf{85} (1982), no. 4, 549--561.
		
	\bibitem{CW} X. Chang, Z.-Q. Wang,
	\emph{Ground state of scalar field equations involving a fractional Laplacian with general nonlinearities},
	Nonlinearity \textbf{26} (2013), no. 2, 479--494.

	\bibitem{CO} Y. Cho, T. Ozawa,
	\emph{Sobolev inequalities with symmetry}, 
	Commun. Contemp. Math. \textbf{11} (2009), no. 3, 355--365.
	
	\bibitem{CT} S. Cingolani, K. Tanaka, 
	\emph{Ground state solutions for the nonlinear Choquard equation with prescribed mass}, 
	to appear on ``Geometric Properties for Parabolic and Elliptic PDE’s'', INdAM Springer Series, Cortona 2019.
	
	\bibitem{DPV} E. Di Nezza, G. Palatucci, E. Valdinoci, 
	\emph{Hitchhiker's guide to the fractional Sobolev spaces}, 
	Bull. Sci. Math. 136 (2012), no. 5, 521--573.
	
	\bibitem{DLQ} L. Dong, D. Liu, W. Qi, L. Wang, H. Zhou, P. Peng, C. Huang,
	\emph{Necklace beams carrying fractional angular momentum in fractional systems with a saturable nonlinearity},
	Commun. Nonlinear Sci. Numer. Simul. \textbf{99} (2021), article ID 105840, pp. 8.

	\bibitem{DZ} S. Duo, Y. Zhang,
	\emph{Computing the ground and first excited states of the fractional Schr\"odinger equation in an infinite potential well}, 
	Commun. Comput. Phys. \textbf{18} (2015), no. 2, 321--350.
	
	\bibitem{FQT} P. Felmer, A. Quaas, J. Tan,
	\emph{Positive solutions of the nonlinear Schr\"odinger equations with the fractional Laplacian},
	Proc. R. Soc. Edinburgh A \textbf{142} (2012), no. 6, 1237--1262.
	
	\bibitem{FJL} J. Fr\"ohlich, B.L.G. Jonsson, E. Lenzmann, 
	\emph{Boson stars as solitary waves}, 
	Comm. Math. Phys. \textbf{274} (2007), no. 1, 1--30.
	
	\bibitem{HIT} J. Hirata, N. Ikoma, K. Tanaka,
	\emph{Nonlinear scalar field equations in $\R^N$: mountain pass and symmetric mountain pass approaches}, 
	Topol. Methods Nonlinear Anal. \textbf{35} (2010), no. 2, 253--276.
	
	\bibitem{HT} J. Hirata, K. Tanaka,
	\emph{Nonlinear scalar field equations with $L^2$ constraint: Mountain Pass and Symmetric Mountain Pass Approaches}, 
	Adv. Nonlinear Stud. \textbf{19} (2019), no. 2, 263--290.
		
	\bibitem{HLRZ} D. Hundertmark, Y.-R. Lee, T. Ried, V. Zharnitsky,
	\emph{Dispersion managed solitons in the presence of saturated nonlinearity}, 
	Physica D \textbf{356-357} (2017), 65--69.
	
	\bibitem{Iko1} N. Ikoma,
	\emph{Existence of solutions of scalar field equations with fractional operator},
	J. Fixed Point Theory Appl. \textbf{19} (2017), no. 2, 649--690.
	
	\bibitem{Iko2} N. Ikoma,
	\emph{Erratum to: Existence of solutions of scalar field equations with fractional operator}, 
	J. Fixed Point Theory Appl. \textbf{19} (2017), no. 2, 1649--1652.
	
	\bibitem{IT} N. Ikoma, K. Tanaka,
	\emph{A note on deformation argument for $L^2$ normalized solutions of nonlinear Schr\"odinger equations and systems}, 
	Adv. Differential Equations 24 (2019), no. 11-12, 609--646.
	
	\bibitem{IP} A.D. Ionescu, F. Pusateri, 
	\emph{Nonlinear fractional Schr\"odinger equations in one dimension},
	J. Funct. Anal. \textbf{266} (2014), no. 1, 139--176.

	\bibitem{Jea0} L. Jeanjean,
	\emph{Existence of solutions with prescribed norm for semilinear elliptic equations},
	Nonlinear Anal. \textbf{28} (1997), no. 10, 1633--1659.
	
	\bibitem{JT1} L. Jeanjean, K. Tanaka,
	\emph{A remark on least energy solutions in $\R^N$},
	Proc. Amer. Math. Soc. \textbf{131} (2003), no. 8, 2399--2408.

	\bibitem{KLS} K. Kirkpatrick, E. Lenzmann, G. Staffilani, 
	\emph{On the continuum limit for discrete NLS with long-range lattice interactions},	
	Comm. Math. Phys. \textbf{317} (2013), no. 3, 563--591.

	\bibitem{KSM} C. Klein, C. Sparber, P. Markowich, 
	\emph{Numerical study of fractional nonlinear Schr\"odinger equations}, 
	Proc. Royal Soc. A \textbf{470} (2014), no. 470, article ID 20140364, pp. 26.

	\bibitem{Las0} N. Laskin, 
	\emph{Fractional quantum mechanics and L\'evy path integrals},
	{Phys. Rev. A} \textbf{268} (2000), no. 4--6, 56--108.

	\bibitem{Lio0} P.-L. Lions, 
	\emph{Sym\'{e}trie et compacit\'e dans les espaces de Sobolev},
	J. Funct. Anal. \textbf{49} (1982), no. 3, 315--334.
	
	\bibitem{Lon0} S. Longhi, 
	\emph{Fractional Schr\"odinger equation in optics} ,
	Optics Lett. \textbf{40} (2015), no. 6, 1117--1120.
		
	\bibitem{MMP} L.A. Maia, E. Montefusco, B. Pellacci,
	\emph{Weakly coupled nonlinear Schr\"odinger systems: the saturation effect},
	Calc. Var. Partial Differential Equations \textbf{46} (2013), no. 1-2, 325--351.

	\bibitem{LMM} P. Li, B. A. Malomed, D. Mihalache,
	\emph{Vortex solitons in fractional nonlinear Schr\"odinger equation with the cubic-quintic nonlinearity},
	Chaos, Solitons and Fractals \textbf{137} (2020), article ID 109783, pp. 10.
	
	\bibitem{Par0} Y. J. Park, 
	\emph{Fractional Gagliardo-Nirenberg inequality},
	J. Chungcheong Math. Soc. \textbf{24} (2011), no. 3, 583--586.

	\bibitem{Rab0} P. H. Rabinowitz,
	``Minimax methods in critical point theory with applications to differential equations'',
	CBMS Reg. Conf. Ser. Math. 65, American Mathematical Society, Providence, 1986.
 
	\bibitem{Shi0} M. Shibata, 
	\emph{Stable standing waves of nonlinear Schr\"odinger equations with a general nonlinear term},
	Manuscripta Math. \textbf{143} (2014), no. 1-2, 221--237.
	
	\bibitem{Stu0} C. Stuart,
	\emph{Bifurcation from the continuous spectrum in the $L^2$-theory of elliptic equations on $\R^n$}, in:
	Recent Methods in Nonlinear Analysis and Applications, Liguori, Naples (1980), 231--300.	
	
	\bibitem{Vaz0} J.L. V\'azquez, 
	\emph{The mathematical theories of diffusion: nonlinear and fractional diffusion}, 
	in ``Nonlocal and Nonlinear Diffusions and Interactions: New Methods and Directions'', Springer, 2017, 205--278.

	\bibitem{ZHSWW} R. Zhang, Z. Han, Y. Shao, Z. Wang, Y. Wang,
	\emph{The numerical study for the ground and excited states of fractional Bose–Einstein condensates},
	Computers and Mathematics with Applications \textbf{78} (2019), 1548---1561.
		
	\bibitem{WCCYZH} Z. Wu, S. Cao, W. Che, F. Yang, X. Zhu, Y. He, 
	\emph{Solitons supported by parity-time-symmetric optical lattices with saturable nonlinearity in fractional Schr\"odinger equation}, 
	Results in Physics \textbf{19} (2020), no. 6, article ID 103381, pp. 6.

	\bibitem{YL} X. Yao, X. Liu,
	\emph{Solitons in the fractional Schr\"odinger equation with parity-time-symmetric lattice potential},
	Photonics Res. \textbf{6} (2018), no. 9, 875--879.
	
\end{thebibliography}
\end{document}